\documentclass[10pt,reqno,a4paper]{amsart}
\usepackage{amssymb}
\usepackage{graphicx}

\newcommand{\om}{\omega}
\newcommand{\omp}{{\omega{\scriptstyle'}}}
\newcommand{\etap}{{\eta{\scriptstyle'}}}

\newcommand{\ri}{\mathrm{i}}
\newcommand{\re}{\mathrm{e}}
\newcommand{\ds}{\displaystyle}
\newcommand{\ded}{{\widehat\eta}}
\newcommand{\wpp}{{\wp{\scriptstyle'}}}

\newcommand{\elambda}{{e_{\scriptscriptstyle\!\lambda}}}
\newcommand{\emu}{{e_{\scriptscriptstyle\!\mu}}}
\newcommand{\sigmalambda}{{\sigma_{\scriptscriptstyle\!\lambda}^{}}}
\newcommand{\sss}{\scriptscriptstyle}

\newcommand{\mfrac}[2]
{\raisebox{0.010em}{\mbox{\footnotesize$\displaystyle\frac{#1}{#2}$}}}

\newcommand{\spin}[1]{\mbox{\footnotesize$\langle#1\rangle$}}

\newcommand{\thetaAB}[2]
{\theta
\raisebox{0.04em}{\mbox{$ \scalebox{0.5}[0.9]{\big[}  \!\!\raisebox{0.1em}
{\scalebox{0.8}{\mbox{\tiny$
\begin{array}{c}#1\\#2\end{array}$}}}
\!\! \scalebox{0.5}[0.9]{\big]}$}}}

\newcommand{\varthetaAB}[2]
{\vartheta
\raisebox{0.04em}{\mbox{$ \scalebox{0.5}[0.9]{\big[}  \!\!\raisebox{0.1em}
{\scalebox{0.8}{\mbox{\tiny$
\begin{array}{c}#1\\#2\end{array}$}}}
\!\! \scalebox{0.5}[0.9]{\big]}$}}}

\begin{document}

\begin{abstract}
We present some new results in theory of classical
$\theta$-functions of Jacobi and $\sigma$-functions of Weierstrass:
ordinary differential equations (dynamical systems) and series
expansions. The paper is basically organized as a stream of new
formulas and constitutes, in re-casted form, some part of author's
own handbook for elliptic/modular functions.
\end{abstract}

\email{brezhnev@mail.ru}
\author{Yu.\,V. Brezhnev}
\title{On functions of Jacobi and Weierstrass (I)}

\maketitle

\section{Introduction}

\noindent Theta-functions of Jacobi and the basis of Weierstrassian
functions $\sigma,\zeta,\wp,\wp'$ occur in numerous theories and
applications. Since their arising, this field was a subject of
intensive studies  and towards the end of  XIX century the state of
the art in this field could be characterized about as follows. If
some formula is not of particular interest then it has already been
obtained. The majority of results were obtained in works of Jacobi
and Weierstrass themselves and  their contemporaries and followers
(Hermite, Kiepert, Hurwitz, Halphen, Frobenius, Fricke, and others).
There is extensive bibliography, including handbooks, on the theory
of elliptic/modular and related functions.  Most thorough works are
four volumes by Tannery \& Molk \cite{tannery}, two volume set by
Halphen with a posthumous edition of the third volume \cite{halphen}
and, of course, {\em Werke\/} by Weierstrass \cite{we} and Jacobi
\cite{jacobi}. A large number of significant examples can be found
in  {\em ``A Course of Modern Analysis''\/} by Whittaker.

In this work we would like to present some  results concerning the
$\theta$-functions of Jacobi and $\sigma$-functions of Weierstrass.
Namely,  power series expansions of the functions and differential
equations to be obeyed by them. In a separate work we shall expound
another way of construction of the Jacobi--Weierstrass theory
different from Jacobi's approach of inversion of a holomorphic
elliptic integral and Weierstrass's construction of elliptic
functions as doubly-periodic and meromorphic ones.

Series expansion of elliptic/modular functions is a widely exploited
apparatus which, due to analytic properties, makes it possible to
get  exact results. It is suffice to mention the function series for
various $\theta$-quotients, the famous number-theoretic $q$-series,
and their consequences like ``Moonshine Conjecture'' and
McKay--Thompson series.

Proposed dynamical systems and their solutions have an independent
interest because solutions of physically significant differential
equations, if these are expressed through Jacobi's $\theta$- or
modular functions, are consequences of the equations presented
below. First of all we should mention the theory of integrable
nonlinear equations and dynamical systems. Differential properties
of the ``modular part'' of Jacobi's functions are more
transcendental and diverse. They gave rise to nice applications over
the last decade known as monopoles and dynamical systems of the
Halphen--Hitchin type \cite{hitchin}, Chazy--Picard--Fuchs equations
considered in works by Harnad, McKay, Ohyama and others. Our
interest in this part of the problem was initiated by a problem of
analytic description of moduli space of algebraic curves and
differential-geometric structures on it. This problem is completely
open, so that the offered ``differential technique'' can be useful
in those situations where some classes of algebraic curves cover
elliptic tori. We shall not touch this topic and restrict ourselves
only to ``book-keeping'' part of the apparatus. It is of interest in
its own right. Applications suggest themselves from the formulas.

In sections 3--10 we shall be keeping the following rule. Unless one
mentioned explicitly or if nothing has been said  about formula then
the result is new\footnote{This means that the formula was not found
in the literature, but does not mean that it is really new one.}. In
some cases, even though a result is obvious, consequences can be
nontrivial. In particular, this concerns the differential equations
on $\theta$-functions with respect to their first argument.
Differential identities sometimes occur in old literature (see for
example \cite{tannery}), but the differential closure is needed. In
most cases the derivation of formulas is self-suggested once a given
formula(s) has been displayed. For this reason we omit proofs and
restrict ourselves to comments to applications or elucidations some
technical details which are important for calculations.

Content of section 2 is standard and presented here to fix notation.
List of bibliography, modern and classical, reduced to a minimum
since, even in a shortened form, it would require mentioning tens
important works.

\section{Definitions and notations}

\subsection{Functions of Jacobi}
Four functions $\theta_{1,2,3,4}^{}$  and their equivalents under
notation with characteristics are defined by the following series:
\begin{center}
\begin{tabular}{l}
$
\begin{array}{l}
\ds\!\!\!-\theta_{11}^{}:\qquad \theta_{1\!}^{}(x|\tau)= -\ri\;
{\sum\limits_{{}^{\scriptscriptstyle k=-\infty}}^{\scriptscriptstyle
\infty}}^{\ds\mathstrut}\! (-1)^k \re^{(k+\frac12)^2\pi\ri\tau}\,
\re^{(2k+1)\pi \ri x}\\
\end{array}
$\\
$
\begin{array}{l}
\ds\phantom{\!\!\!-}\theta_{10}^{}:\qquad \theta_{2\!}^{}(x|\tau)=
\phantom{-\ri\;} {\sum\limits_{{}^{\scriptscriptstyle
k=-\infty}}^{\scriptscriptstyle \infty}}^{\ds\mathstrut}
\re^{(k+\frac12)^2\pi\ri\tau}\,\re^{(2k+1)\pi \ri x}
\end{array}
$\\
$
\begin{array}{l}
\ds\phantom{\!\!\!-}\theta_{00}^{}:\qquad \theta_{3\!}^{}(x|\tau)=
\phantom{-\ri\;}{\sum\limits_{\scriptscriptstyle
k=-\infty}^{\scriptscriptstyle \infty}}^{\ds\mathstrut}
\re^{k^2\pi\ri\tau}\,\re^{2k\pi \ri x}
\end{array}
$\\
$
\begin{array}{l}
\ds\phantom{\!\!\!-}\theta_{01}^{}:\qquad \theta_{4\!}^{}(x|\tau)=
\phantom{-\ri\;}{\sum\limits_{{}^{\scriptscriptstyle
k=-\infty}}^{\scriptscriptstyle \infty}}^{\ds\mathstrut}\! (-1)^k
\re^{k^2\pi\ri\tau}\, \re^{2k\pi \ri x}
\end{array}
$\\
\end{tabular}.
\end{center}
We shall use the designation: $\theta_k\equiv \theta_k(x|\tau)$.
Values of $\theta$-functions under $x=0$ are called
$\vartheta$-constants:
$\vartheta_k\equiv\vartheta_k(\tau)=\theta_k(0|\tau)$. For handy
formatting formulas we shall use twofold notation for
$\theta$-functions with characteristics:
$$
\theta \raisebox{0.04em}{\mbox{$ \scalebox{0.5}[0.9]{\big[}
\!\!\raisebox{0.1em}
{\scalebox{0.8}{\mbox{\tiny$\begin{array}{c}\alpha\\\beta
\end{array}$}}}
\!\! \scalebox{0.5}[0.9]{\big]}$}}(x|\tau)=
\theta_{\alpha\beta}(x|\tau)= {\sum\limits_{{}^{\scriptscriptstyle
k=-\infty}}^{\scriptscriptstyle \infty}}^ {\ds\mathstrut}
\re^{\pi\ri \left(\!k+\frac\alpha2\!\right)^{\!2}\tau+
2\pi\ri\left(\!k+\frac\alpha2\!\right)\!
\left(\!x+\frac\beta2\!\right)\! }_{\mathstrut}\,.
$$
We consider arbitrary integral characteristics and therefore the
functions $\theta_{\alpha\beta}$ are always equal to
$\pm\theta_{1,2,3,4}$. Shifts by $\frac12$-periods lead to shifts of
characteristics:
$$
\thetaAB{\alpha}{\beta} \!\left({\textstyle x+\frac n2+\frac
m2\,\tau}\big|\tau\right)=
\thetaAB{\alpha\mbox{\tiny+}m}{\beta\mbox{\tiny+}n}(x|\tau)
\!\cdot\!\re^{\mbox{\tiny--}\pi\ri
\,m\left(x+\frac{\beta+n}{2}+\frac m4\tau\right)} _{\mathstrut}\ ,
\qquad\qquad(n,\,m)=0,\pm1,\pm2,\ldots
$$
Twofold shifts by $\frac12$-periods yield the law of transformation
of $\theta$-function into itself:
$$
\theta_{\alpha\beta} (x+n+m\,\tau|\tau)= \theta_{\alpha\beta}
(x|\tau)\!\cdot\! \re^{\mbox{\tiny--}\pi\ri\,(m^2\tau+2mx)}\,
(-1)^{n\alpha-m\beta}\,,\qquad\qquad(n,\,m)=0,\pm1,\pm2,\ldots.
$$
Value of any  $\theta$-function at any $\frac12$-period is a certain
$\vartheta$-constant with an exponential multiplier:
$$
\thetaAB{\alpha}{\beta} \!\left({\textstyle\frac n2+\frac
m2\,\tau}\big|\tau\right)=
\varthetaAB{\alpha\mbox{\tiny+}m}{\beta\mbox{\tiny+}n}(\tau)
\!\cdot\!\re^{\mbox{\tiny--}\pi\ri m\left(\frac{\beta+n}{2}+\frac
m4\tau\right)} _{\mathstrut}\ ,
\qquad\qquad\qquad(n,\,m)=0,\pm1,\pm2,\ldots
$$

\subsection{Functions of Weierstrass}
Weierstrassian notation $\sigma(x|\om,\omp)$,
$\wp(x|\om,\omp),\,\ldots$ and $\sigma(x;g_2^{},g_3^{})$,
$\wp(x;g_2^{},g_3^{}),\, \ldots$ is the commonly accepted. Due to
known relations of homogeneity for functions
$\sigma,\,\zeta,\,\wp,\,\wpp$, the two parameters $(\om,\omp)$ or
$(g_2^{},g_3^{})$ can be replaced by one $\tau=\frac{\omp}{\om}$
which is called modulus of elliptic curves. We shall designate that
functions as follows
$$
\sigma(x|\tau)\equiv\sigma(x|1,\tau),\quad
\zeta(x|\tau)\equiv\zeta(x|1,\tau)\quad
\wp(x|\tau)\equiv\wp(x|1,\tau)\quad
\wpp(x|\tau)\equiv\wpp(x|1,\tau)\,.
$$
Weierstrassian invariants $g_2^{},\,g_3^{}$ and $\Delta =
g_2^3-27g_3^2$ are functions of periods or modulus and defined by
known formulas of Weierstrass--Eisenstein. These series are entirely
unsuited for numeric computations. Hurwitz, in his dissertation
(1881), found a nice transition to function series. Such series are
most effective for computations:
$$
\begin{array}{l}
\ds g_2^{}(\tau) =20\,\pi^4 \left\{\frac{1}{240}
+\sum_{{}^{\scriptscriptstyle k=1}}^{\scriptscriptstyle \infty}
\frac{k^3\,\re^{2k\pi\ri\tau}}{1-\re^{2k\pi\ri\tau}_{\mathstrut}}
\right\},\qquad  g_3^{}(\tau) = \mfrac73\,\pi^6 \left\{\frac{1}{504}
- \sum_{{}^{\scriptscriptstyle k=1}}^{\scriptscriptstyle \infty}
\frac{k^5\,\re^{2k\pi\ri\tau}_{}}{1-\re^{2k\pi\ri\tau}_{\mathstrut}}
\right\}\,.
\end{array}
$$
$\eta$-function of Weierstrass is defined by the formula
$\eta(\tau)=\zeta(1|1,\tau)$ and the following formula is used for
computations:
$$
\eta(\tau)=2\pi^2 \left\{ \frac{1}{24}- \sum_{{}^{\scriptscriptstyle
k=1}}^{\scriptscriptstyle \infty}
\frac{\re^{2k\pi\ri\tau}}{\big(1-\re^{2k\pi\ri\tau}_{\mathstrut}\big)^2}
\right\}, \qquad
\ds\eta\!\left(\mfrac{a\,\tau+b}{c\,\tau+d}\right)=\ds
(c\tau+d)^2\,\eta(\tau)-\frac{\pi\ri}{2}\,c\,(c\,\tau+d)\,,
$$
where the numbers $(a,b,c,d)$ are integer and $ad-bc=1$. Modular
transformations are used not only in theory but also in practice
since the value of modulus strongly affects the convergence of the
function series. Moving $\tau$  into a fundamental domain of modular
group (the process is easily automatized), one obtains the values of
$\tau$ with the minimally allowable imaginary part $\Im(\tau)=
\frac{\sqrt{3}}{2}$. At such  ``the most worst'' point the series
converge very fast.

Three functions of Weierstrass $\sigmalambda$ are defined by the
following expressions:
$$
\begin{array}{l}
\ds \sigmalambda\!(x|\om,\omp)=
\phantom{-}\re^{\eta(\om,\om{\sss'})\,\frac{x^2}{2\om}}_{\mathstrut}\!\cdot\!
\frac{\theta_{\sss\lambda+1}\!
\big(\frac{x}{2\om}\big|\frac{\om{\sss'}}{\om}\big)}
{\vartheta_{\sss\lambda+1}\! \big(\frac{\om{\sss'}}{\om}\big)}\;,
\qquad\mbox{where \ \ $\lambda=1,2,3$}\,.
\end{array}
$$
The $\sigma$-function of Weierstrass, as a function of
$(x,\,g_2^{},\,g_3^{})$, satisfies the linear differential equations
obtained by Weierstrass:
\begin{equation}\label{s12}\left\{
\begin{array}{r}
\ds
x\,\frac{\partial\sigma}{\partial x}-
4\,g_2^{}\,\frac{\partial\sigma}{\partial g_2^{}}-
6\,g_3^{}\,\frac{\partial\sigma}{\partial g_3^{}}-\sigma=0\;\,\\ \\
\ds \frac{\partial^2\sigma}{\partial x^2}-
12\,g_3^{}\,\frac{\partial\sigma}{\partial g_2^{}}-
\mfrac23\,g_2^2\,\frac{\partial\sigma}{\partial g_3^{}}
+\mfrac{1}{12}\,g_2^{}\,x^2\,\sigma=0
\end{array}\right..
\end{equation}
These equations make it possible to write down the recursive
relation for coefficients $C_k(g_2^{},g_3^{})$ of the power series
of the function $\sigma$:
\begin{equation}\label{sigma}
\sigma(x;g_2^{},g_3^{})= \displaystyle
C_0\,x+C_1\,\mfrac{x^3}{3!}+\cdots=
x-\mfrac{g_2^{}}{240}\,x^5-\mfrac{g_3^{}}{840}\,x^7 +\cdots, \qquad
\big( C_0=1,\; C_1=0\big)\,.
\end{equation}
Two such recurrences are known. One is due to Halphen:
\begin{equation}\label{halphen}
C_k=-\mbox{\large$\widehat{\boldsymbol{\mathfrak D}}$}\,C_{k-1}
-\mfrac16\,(k-1)(2\,k-1)\,g_2^{}\,C_{k-2}\;,\quad \mbox{where}\quad
 \mbox{\large$\widehat{{\mathfrak D}}$}=
-12\,g_3^{}\,\frac{\partial}{\partial g_2^{}}-
\mfrac23\,g_2^2\,\frac{\partial}{\partial g_3^{}}\,.
\end{equation}
The second was obtained by Weierstrass:
$$
\sigma(x;g_2^{},g_3^{})=\ds
\mbox{\large$\ds\sum_{\sss m,n=0}^{\sss\infty}$}
\,A_{m,\,n}
\left(\!{\mfrac{g_2^{}}{2}}\!\right)^{\!m}
\big( 2g_3^{}\big)^n {\mfrac{x^{4m+6n+1}}
{(4m+6n+1)!}},
\qquad
\left\{
\begin{array}{rcl}
A_{\sss 0,0}\!\!\!&=&\!\!\!1\\
A_{m,\,n}\!\!\!&=&\!\!\!0\quad
\big(m<0,\; n<0\big)
\end{array}
\right..
$$
$$
A_{m,\,n}=
\mfrac{16}{3}\,(n+1)\,A_{m-2,\,n+1}
+3\,(m+1)\,A_{m+1,\,n-1}
-\mfrac13\, (2m+3n-1)(4m+6\,n-1)\,A_{m-1,\,n}.
$$
Recently, in connection with theory of Kleinian
$\boldsymbol{\sigma}$-functions, yet another recurrence was obtained
\cite{eilbeck}. Among all the recurrences, Weierstrassian one is the
least expendable. There are only multiplications of integers in it
and its comparison efficiency rapidly grows under the growing of
order of the expansions. Weierstrass proves separately that the
numbers $A_{m,n}$ are integers.

\subsection{Function of Dedekind}
Due to coincidence the standard designations for Weierstrassian
function $\eta(\tau)$ and Dedekind's one, we shall use, for the last
one, the sign $\ded(\tau)$:
$$
\ded(\tau)= \re^{\frac{\pi\ri}{12}\tau}\,
{\prod\limits_{\scriptscriptstyle k=1}^{\scriptscriptstyle \infty}}^
{\ds\mathstrut}_{\ds\mathstrut} \big(1-\re^{2k\pi\ri\tau}\big)=
\re^{\frac{\pi\ri}{12}\tau}_{} {\sum\limits_{{}^{\scriptscriptstyle
k=-\infty}}^{\scriptscriptstyle \infty}}^ {\ds\mathstrut}
\!(-1)^k\,\re^{(3k^2+k)\pi\ri\tau}\,.\qquad (\mbox{Euler (1748)})
$$
There is a differential relation between the two functions. It is
expressed by the formula $\mfrac{1}{\ded}\,\mfrac{d
\ded}{d\tau}=\mfrac{\ri}{\pi}\,\eta$.

\section{Differentiation of functions
$\sigma,\zeta,\wp,\wpp$}

\noindent Differentiations  of Weierstrassian functions with respect
to invariants are known \cite{halphen,tannery}. Rules of
transformations between differentiations with respect to
$(g_2^{},g_3^{})$ and $(\om,\omp)$ are also known
(Frobenius--Stickelberger (1882)) \cite{halphen}. Applying them, and
assuming $\Im\big(\frac{\omp}{\om}\big)>0$, one obtains
$$
\left\{
\begin{array}{rcl}
\displaystyle  \frac{\partial\sigma}{\partial\om}\!\!\!&=&\!\!\!
\displaystyle -\frac{\ri}{\pi} \left\{
\omp\Big(\wp-\zeta^2-\mfrac{1}{12}\,g_2^{}\,x^2 \Big)
+2\,\etap\,\big(x\,\zeta-1\big)\right\}\sigma\\ \\
\displaystyle  \frac{\partial\sigma}{\partial\omp}\!\!\!&=&\!\!\!
\displaystyle \phantom{-}\frac{\ri}{\pi} \left\{
\om\Big(\wp-\zeta^2-\mfrac{1}{12}\,g_2^{}\,x^2 \Big)
+2\,\eta\,\big(x\,\zeta-1\big)\right\}\sigma
\end{array}
\right.,
$$
\medskip
$$
\left\{
\begin{array}{rcl}
\displaystyle  \frac{\partial\zeta}{\partial\om}\!\!\!&=&\!\!\!
\displaystyle -\frac{\ri}{\pi}\, \left\{
\omp\Big(\wpp+2\,\zeta\,\wp-\mfrac16\,g_2^{}\,x\Big)
+2\,\etap\big(\zeta-x\,\wp\big)\right\}\\\\
\displaystyle  \frac{\partial\zeta}{\partial\omp}\!\!\!&=&\!\!\!
\phantom{-}\displaystyle \frac{\ri}{\pi}\, \left\{
\om\,\Big(\wpp+2\,\zeta\,\wp-\mfrac16\,g_2^{}\,x\Big)
+2\,\eta\,\big(\zeta-x\,\wp\big)\right\}
\end{array}
\right.,
$$
\medskip
$$
\;\,\left\{
\begin{array}{rcl}
\displaystyle  \frac{\partial\wp}{\partial\om}\!\!\!&=&\!\!\!
\phantom{-}\displaystyle \frac{2\,\ri}{\pi}\, \left\{
\omp\Big(2\,\wp^2+\zeta\,\wpp-\mfrac13\,g_2^{}\Big)
-\etap\big(2\,\wp+x\,\wpp\big)\right\}\\ \\
\displaystyle  \frac{\partial\wp}{\partial\omp}\!\!\!&=&\!\!\!
\displaystyle -\frac{2\,\ri}{\pi}\, \left\{
\om\,\Big(2\,\wp^2+\zeta\,\wpp-\mfrac13\,g_2^{}\Big)
-\eta\,\big(2\,\wp+x\,\wpp\big)\right\}
\end{array}
\right.,
$$
$$
\left\{
\begin{array}{rcl}
\displaystyle  \frac{\partial\wpp}{\partial\om}\!\!\!&=&\!\!\!
\phantom{-}\displaystyle \frac{\ri}{\pi}\, \Big\{
\omp\big(6\,\wp\,\wpp+12\,\zeta\,\wp^2-g_2^{}\,\zeta\big)
-\etap\big(6\,\wpp+12\,x\,\wp^2-g_2^{}\,x\big)\Big\}\\ \\
\displaystyle  \frac{\partial\wpp}{\partial\omp}\!\!\!&=&\!\!\!
\displaystyle -\frac{\ri}{\pi}\, \Big\{
\om\,\big(6\,\wp\,\wpp+12\,\zeta\,\wp^2-g_2^{}\,\zeta\big)
-\eta\,\big(6\,\wpp+12\,x\,\wp^2-g_2^{}\,x\big)\Big\}
\end{array}
\right.\,.
$$
Setting here $(\om=1, \,\omp=\tau)$ one obtains the dynamical system
with a parameter $x$:
$$ \left\{
\begin{array}{rcl}
\ds\frac{\partial\sigma}{\partial\tau}\!\!\!&=&\!\!\! \displaystyle
\phantom{-2}\frac{\ri}{\pi} \left\{ \wp-\zeta^2
+2\,\eta\,(x\,\zeta-1)-\mfrac{1}{12}\,g_2^{}\,x^2\right\}\sigma\\\\
\ds\frac{\partial\zeta}{\partial\tau}\!\!\!&=&\!\!\!
\phantom{-2}\displaystyle \frac{\ri}{\pi} \left\{ \wpp+
2\,\eta\zeta+2\,\wp\,(\zeta-x\,\eta)
-\mfrac16\,g_2^{}\,x\right\}\\\\
\ds\frac{\partial\wp}{\partial\tau}\!\!\!&=&\!\!\! \ds
-2\frac{\ri}{\pi} \left\{ 2\,\wp^2+\wpp(\zeta-x\,\eta)
-2\,\eta\,\wp-\mfrac13\,g_2^{}\right\}\\\\
\ds\frac{\partial\wpp}{\partial\tau}\!\!\!&=&\!\!\! \ds
-6\frac{\ri}{\pi} \left\{
\wpp(\wp-\eta)+\Big(2\,\wp^2-\mfrac16\,g_2^{}\Big)(\zeta-x\,\eta)
\right\}
\end{array}
\right.\,.
$$
The function $\sigma(x|\tau)$ satisfies one differential equation
$$
\frac{\partial^2\sigma}{\partial x^2}-2\,x\,\eta\,
\frac{\partial\sigma}{\partial x}-\pi\,\ri\,
\frac{\partial\sigma}{\partial \tau}+
\left(2\,\eta+\mfrac{1}{12}\,g_2^{}\,x^2 \right)\sigma=0\,,
$$
which is an analogue of heat equation for the function
$\theta(x|\tau)$. All the differential equations in variables
$x,\om,\omp,\tau$  contain no the quantity $g_3^{}$. Depending on
representation $(g_2^{},g_3^{})$, $(\om,\omp)$ or $\tau$, it is
either algebraic integral of these equations
$g_3^{}=4\,\wp^3-g_2^{}\,\wp-\wpp^2$ or determines a fix algebraic
relation between the variables. Each of the functions
$\sigma,\zeta,\wp,\wpp(x|\tau)$ satisfies ordinary differential
equation of 4-th order  with variable coefficients
$g_{2,3}^{}(\tau),\, \eta(\tau)$ in the variable $\tau$. These
equations have defied simplifications. By virtue of linear
connection between function $\sigma$ and $\theta$, the last one also
satisfies ordinary differential equations.

\section{Power series for $\sigma$-functions of Weierstrass}

\noindent Weierstrassian recurrence has the following graphical
interpretation:
\medskip

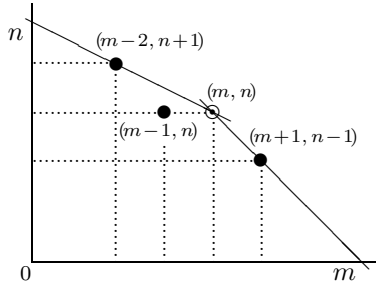
\centerline{\unitlength=0.8mm
\begin{picture}(60,60)
\put(0,15){\line(1,0){57}} \put(0,15){\line(0,1){43}}
\multiput(29.7,15)(0,1){24}{\scriptsize .}
\multiput(0,39.75)(1,0){21}{\scriptsize .}
\multiput(23.3,39.75)(1,0){6}{\scriptsize
.}\put(22,40){\circle{0.1}} \put(30,40){\circle{2}}
\put(30,40){\circle*{0.9}} \put(29,42.5){{\scriptsize$(\!m,n\!)$}}
\put(22,40){\circle*{2}}
\put(14.3,36){{\scriptsize$(\!m\!-\!1,n\!)$}}
\multiput(21.5,15)(0,1){20}{\scriptsize .} \put(14,48){\circle*{2}}
\put(10.5,51){{\scriptsize$(\!m\!-\!2,n\!+\!1\!)$}}
\multiput(0,48)(1,0){13}{\scriptsize .}
\multiput(13.5,15)(0,1){33}{\scriptsize .} \put(14,48){\circle{0.1}}
\put(38,32){\circle{0.1}}\put(38,32){\circle*{2}}
\put(36,35){{\scriptsize$(\!m\!+\!1,n\!-\!1\!)$}}
\multiput(37.7,15)(0,1){18}{\scriptsize .}
\multiput(0,31.7)(1,0){37}{\scriptsize .} \put(50,12){$m$}
\put(-4,52){$n$} \put(-2,12){\footnotesize $0$}
\put(28,42){\line(1,-1){28}} \put(-1,55.5){\line(2,-1){34}}
\end{picture}
} \vspace{-0.5cm}

\noindent The explicitly grouped  $\sigma$-series (\ref{sigma}) has
the form ($[n]$ denotes an integral part of the number)
\begin{equation}\label{wei}
\sigma(x;g_2^{},g_3^{})= \mbox{\large$\ds\sum_{\scriptstyle
k=0}^\infty$}\, \left\{ \sum_{\sss\nu=[k/3]}^{\sss k/2}
\!2^{2k-5\nu} A_{3\nu-k,\,k-2\nu}\,\cdot g_2^{3\nu-k}\, g_3^{k-2\nu}
\right\} \frac{x^{2k+1}}{(2k+1)!}\ .
\end{equation}

Denote $\elambda\equiv \wp(\om_{\!\sss\lambda}^{}|\om,\omp)$. Then
the functions $\sigmalambda$ satisfy the equations of Halphen:
$$
\left\{
\begin{array}{rcl}
\ds
x\,\frac{\partial\sigmalambda}{\partial x}-
2\,\elambda\,\frac{\partial\sigmalambda}{\partial\elambda}-
4\,g_2^{}\,\frac{\partial\sigmalambda}{\partial g_2^{}}\!\!\!&=&\!\!\!0\\ \\
\ds \frac{\partial^2\sigmalambda}{\partial x^2}-
\left(4\,\elambda^{\!\!2}-\mfrac23\,g_2^{}\right)\!
\frac{\partial\sigmalambda}{\partial \elambda}-
12\,\big(4\,\elambda^{\!\!3}-g_2^{}\,\elambda\big)
\frac{\partial\sigmalambda}{\partial g_2^{}}+
\left(\elambda+\mfrac{1}{12}\,g_2^{}\,x^2\right)\!
\sigmalambda\!\!\!&=&\!\!\!0
\end{array}
\right..
$$
Herefrom one obtains analogue of the recurrence (\ref{wei}):
$$
\begin{array}{rl}
\ds
\sigmalambda\!(x;\elambda,g_2^{})=&\ds
\!\!
\mbox{\large$\ds\sum_{\scriptstyle k=0}^\infty$}
\left\{
\sum_{\sss \nu=0}^{\sss k/2}\,2^{-\nu}\,\mathfrak{B}_{k-2\nu,\,\nu}\cdot
\elambda^{\!\!k-2\nu}g_2^\nu
\right\}\frac{x^{2k}}{(2k)!}\, ,\\ \\
\ds
\mathfrak{B}_{m,n}=&\!\!24\,(n+1)\,\mathfrak{B}_{m-3,\,n+1}+
(4\,m-12\,n-5)\,\mathfrak{B}_{m-1,\,n}-\\ \\
&\!\!\ds
-\mfrac43\,(m+1)\,\mathfrak{B}_{m+1,\,n-1}-
\mfrac13(m+2\,n-1)(2\,m+4\,n-3)\,\mathfrak{B}_{m,\,n-1}\, .
\end{array}
$$
All the coefficients $\mathfrak{B}_{m,\,n}$ are integrals, and
$\mathfrak{B}_{\sss 0,0}=1$ and $\mathfrak{B}_{m,\,n}=0$ under
$(m,n)<0$. Weierstrass wrote out recurrences for the functions
$S_\lambda^{}=\re^{\frac12\elambda x^2}_{}\sigmalambda$ in
representation $\big(\elambda,\,\varepsilon_{\sss\lambda}=
3\elambda^{\!\!2}-\frac14g_2^{}\big)$. One can build a universal
series for the functions $\sigma,\sigmalambda$. All the functions
satisfy the differential equations
\begin{equation}\label{epsilon}
\left\{
\begin{array}{rcl}
\ds
x\,\frac{\partial\Xi}{\partial x}-
2\,\elambda\,\frac{\partial\Xi}{\partial\elambda}-
4\,g_2^{}\,\frac{\partial\Xi}{\partial g_2^{}}
-(1-\varepsilon)\,\Xi\!\!\!&=&\!\!\!0\\ \\
\ds \frac{\partial^2\Xi}{\partial x^2}-
\left(4\,\elambda^{\!\!2}-\mfrac23\,g_2^{}\right)\!
\frac{\partial\Xi}{\partial \elambda}-
12\,\big(4\,\elambda^{\!\!3}-g_2^{}\,\elambda\big)
\frac{\partial\Xi}{\partial g_2^{}}+
\left(\varepsilon\,\elambda+\mfrac{1}{12}\,g_2^{}\,x^2\right)\!
\Xi\!\!\!&=&\!\!\!0
\end{array}
\right.,
\end{equation}
where  case $\Xi=\sigmalambda$ corresponds to $\varepsilon=1$, and
$\Xi=\sigma$ corresponds to $\varepsilon=0$  and arbitrary
$\elambda$. The quantity $\varepsilon$ is a parity of
$\theta$-characteristic. The universal series and integer recurrence
acquire the following form:
\begin{equation}\label{universal}
\begin{array}{rl}
\ds
\Xi(x;\elambda,g_2^{})=&\ds
\!\!\mbox{\large$\ds\sum_{k=0}^{\infty}$}
\left\{
\sum_{\sss\nu=0}^{\sss k/2}\,2^{-\nu}\,
\mathfrak{B}_{k-2\nu,\,\nu}^{{\sss(}\varepsilon{\sss)}}\cdot
\elambda^{\!\!k-2\nu}g_2^\nu
\right\}\frac{x^{2k+1-\varepsilon}}
{(2k+1-\varepsilon)!}\, ,\\ \\
\ds
\mathfrak{B}_{m,n}^{{\sss(}\varepsilon{\sss)}}=&
\!\!24\,(n+1)\,\mathfrak{B}_{m-3,\,n+1}^{{\sss(}\varepsilon{\sss)}}+
(4\,m-12\,n-4-\varepsilon)\,\mathfrak{B}_{m-1,\,n}^
{{\sss(}\varepsilon{\sss)}}-\\ \\
&\!\!\ds
-\mfrac43\,(m+1)\,\mathfrak{B}_{m+1,\,n-1}^{{\sss(}\varepsilon{\sss)}}-
\mfrac13 (m+2\,n-1)(2\,m+4\,n-1-2\,\varepsilon)\,\mathfrak{B}_{m,\,n-1}
^{{\sss(}\varepsilon{\sss)}}\, .
\end{array}
\end{equation}
Transition between pairs $(g_2^{},g_3^{})$, $(\elambda,g_2^{})$ and
$(\elambda,\emu)$ is one-to-one therefore the universal recurrence
can be written for any of these representations. In the last case it
is symmetrical. The representation (\ref{universal}) is a five-term
one, in contrast to four-term recurrence of Weierstrass $A_{m,n}$,
but the representation $(g_2^{},g_3^{})$ is not possible for the
functions $\sigmalambda$.

\section{Power series for $\theta$-functions of Jacobi}

\noindent Jacobi made an attempt to obtain such series as early as
before appearance of his {\em Fundamenta Nova\/} \cite[{\bf I}:
259--260]{jacobi}. These series are the series with polynomial
coefficients in $\eta(\tau)$ and $\vartheta(\tau)$-constants. This
fact follows from the formulas
$$
\theta_1^{}\!(x|\tau)=\pi\,\ded^3(\tau)
\!\cdot\!\re^{-2\eta(\tau)x^2}_{\mathstrut}\sigma(2x|\tau)\;,
\qquad\qquad \theta_{\sss
\lambda}\!(x|\tau)=\vartheta_{\sss\!\lambda}^{}\!(\tau)
\!\cdot\!\re^{-2\eta(\tau)x^2}_{\mathstrut} \sigma_{\sss
\!\lambda-1}^{}\!(2x|\tau)\,.
$$
$\vartheta$-constants make it possible to rewrite Weierstrassian
representations $g_2^{},g_3^{},\elambda$ into $\vartheta$-constant
ones.  Formulas for branch-points $e_k^{}$ through the
$\vartheta$-constants are well known. In turn, $\vartheta$-constants
are connected via Jacobi identity
$\vartheta_3^4=\vartheta_2^4+\vartheta_4^4$. All this enables one to
convert the representations choosing  arbitrary pair. We shall speak
{\em $(\alpha,\beta)$-representation\/}, if formulas are written
through the constants
$(\vartheta_{\sss\alpha0},\vartheta_{\sss0\beta})$ under
$(\alpha,\beta)\ne(0,0)$. Operator  (\ref{halphen}) has
Weierstrassian representation
$$
\begin{array}{l}
\ds {\mbox{\large$\widehat{\boldsymbol{\mathfrak D}}$}}=
-12\,\big(4\,\elambda^{\!\!3}-g_2^{}\,\elambda\big)\frac{\partial}{\partial
g_2^{}}- \left(4\,\elambda^{\!\!2}-\mfrac23 g_2^{}\right)\,
\frac{\partial}{\partial\elambda}=\\ \\
\ds \phantom{\mbox{\large$\widehat{\boldsymbol{\mathfrak D}}$}}=
\mfrac43\big(2\,\emu^{\!\!2}+2\,\emu\elambda-\elambda^{\!\!2}\big)
\frac{\partial}{\partial \elambda}+
\mfrac43\big(2\,\elambda^{\!\!2}+2\,\elambda\emu-\emu^{\!\!2}\big)
\frac{\partial}{\partial\emu}\,,
\end{array}
$$
and also $\vartheta$-constant one. For example
$(\vartheta_2^{},\vartheta_4^{})$-representation
\begin{equation}\label{tmp}
\mbox{\large$\widehat{\boldsymbol{\mathfrak D}}$}= \frac{\pi^2}{3}
\left\{\vartheta_2^8+2\,\vartheta_2^4\,\vartheta_4^4\right\}
\frac{\partial}{\partial\vartheta_2^4}- \frac{\pi^2}{3}
\left\{\vartheta_4^8+2\,\vartheta_2^4\,\vartheta_4^4\right\}
\frac{\partial}{\partial\vartheta_4^4}\;.
\end{equation}
Let us denote $\spin{\alpha}=(-1)^\alpha$. General
$(\alpha,\beta)$-representation of Halphen's operator has the form
$$
\mbox{\large$\widehat{\boldsymbol{\mathfrak D}}$}=
\frac{\pi^2}{3}\left(\spin{\alpha}\,
\vartheta_{\sss0\beta}{}^{\!\!\!\!\!\!8}\;+ 2\,\spin{\beta}\,
\vartheta_{\sss\alpha0}{}^{\!\!\!\!\!\!4}\;\;
\vartheta_{\sss0\beta}{}^{\!\!\!\!\!\!4}\;\right)\!
\frac{\partial}{\partial\vartheta_{\sss0\beta}{}^{\!\!\!\!\!\!4}}
-\frac{\pi^2}{3}\left(\spin{\beta}\,
\vartheta_{\sss\alpha0}{}^{\!\!\!\!\!\!8}\;+ 2\,\spin{\alpha}\,
\vartheta_{\sss\alpha0}{}^{\!\!\!\!\!\!4}\;\;
\vartheta_{\sss0\beta}{}^{\!\!\!\!\!\!4}\; \right)\!
\frac{\partial}{\partial\vartheta_{\sss\alpha0}{}^{\!\!\!\!\!\!4}}\,.
$$
Let $\varepsilon$ is defined by a parity of characteristic of the
function $\theta_{\sss\alpha\beta}(x|\tau)$:
$$
\varepsilon=\frac{\spin{\alpha\beta}+1}{2}\,, \qquad\qquad
\left\{\!\!\!
\begin{array}{l}
\varepsilon=0\quad \mbox{under \ } \theta_{\sss\!\alpha\beta}=\pm\,\theta_1^{}\\
\varepsilon=1\quad \mbox{under \ } \theta_{\sss\!\alpha\beta}=\pm\,
\theta_{2,3,4}^{\ds\mathstrut}
\end{array}\!\!\!
\right\}\;.
$$
The equation (\ref{epsilon}) for the function
$\Xi=(\sigma,\sigmalambda)$ has the form
$$
\frac{\partial^2\,\Xi}{\partial x^2}+
\mbox{\large$\widehat{\boldsymbol{\mathfrak D}}$}\, \Xi
+\left\{\varepsilon\,\elambda(\vartheta)+\mfrac{\pi^4}{12^2}
\!\left[\vartheta_2^8+\vartheta_2^4\,\vartheta_4^4+ \vartheta_4^8
\right]\!x^2\right\}\Xi=0\;,
$$
where $\mbox{\large$\widehat{\boldsymbol{\mathfrak D}}$}$ is
determined by formula (\ref{tmp}).

\subsection{Function $\theta_1^{}\!(x|\tau)$} Expansion of function
\begin{equation}\label{theta1}
\theta_1^{}\!(x|\tau)=\sum_{k=0}^{\infty}\,C_k(\tau)
\!\cdot\!x^{2k+1}=2\,\pi\,\ded^3 \left\{
x-2\,\eta\!\cdot\!x^3+\Big(2\,\eta^2- \mfrac{\pi^4}{180}
\big(\vartheta_2^8+\vartheta_2^4\,\vartheta_4^4+\vartheta_4^8\big)\Big)
\!\cdot\!x^5+\cdots \right\}
\end{equation}
has the following analytic representation:
\begin{equation}\label{t1}
\begin{array}{l}
\ds \theta_1^{}\!(x|\tau)=\phantom{\ded^3}2\,\pi\,\,
\mbox{\large$\ds\sum_{\scriptstyle k=0}^\infty$}
\,\frac{(4\pi\ri)^k}{(2k+1)!}\,\frac{d^k\,\ded^3}{d\tau^k}
\!\cdot\!x^{2k+1}=\\ \\
\ds \phantom{\!\theta_1^{}(x|\tau)}= 2\,\pi\,\ded^3\,\,
\mbox{\large$\ds\sum_{\scriptstyle k=0}^\infty$} \,(-2)^k\left\{
\sum_{\sss\nu=0}^{\sss k}
\mfrac{(-6)^{-\nu}\,\pi^{2\nu}}{(k-\nu)!\,(2\nu+1)!} \!\cdot\!
\eta^{k-\nu}\,\mathcal{N}_\nu(\vartheta) \right\} x^{2k+1}\;,
\end{array}
\end{equation}
where polynomials  $\mathcal{N}_\nu(\vartheta)$, depending on
combinations of $\vartheta$-constants, have the form
$$
\mathcal{N}_\nu(\vartheta)= \mbox{\large$\ds\sum_{s=0}^{\nu}$}\;
\mbox{\footnotesize$\textstyle\left\{\!\!\!
\begin{array}{r}
\mbox{\small$\mathfrak{G}_{\nu-s,\,s}$}\cdot\vartheta_4^{4s\mathstrut}\,
\vartheta_2^{4{\sss(}\nu-s{\sss)}}\\
(-1)^s\,\mbox{\small$\mathfrak{G}_{\nu-s,\,s}$}\cdot
\vartheta_3^{4s\ds\mathstrut}\,
\vartheta_4^{4{\sss(}\nu-s{\sss)}\ds\mathstrut}\\
(-1)^s\,\mbox{\small$\mathfrak{G}_{s,\,\nu-s}$}\cdot\vartheta_3^{4s\ds\mathstrut}\,
\vartheta_2^{4{\sss(}\nu-s{\sss)}\ds\mathstrut}
\end{array}\!\!\!
\right\}$}.
$$
Integer recurrence $\mathfrak{G}_{m,n}$ $(\mathfrak{G}_{\sss
0,\,0}=1)$ looks as follows:
$$
\begin{array}{l}
\ds \mathfrak{G}_{m,\,n}=4\,(n-2\,m-1)\,\mathfrak{G}_{m,\,n-1}-
4\,(m-2\,n-1)\,\mathfrak{G}_{m-1,\,n}-
\\ \\
\ds\phantom{\mathfrak{G}_{m,\,n}=}
-2\,(m+n-1)(2\,m+2\,n-1)\big(\mathfrak{G}_{m-2,\,n}+
\mathfrak{G}_{m-1,\,n-1}+\mathfrak{G}_{m,\,n-2}\big)
\end{array}.
$$
The recurrence $\mathfrak{G}_{m,n}$ is antisymmetric:
$\mathfrak{G}_{m,n}=(-1)^{m+n}\,\mathfrak{G}_{n,m}$. One may write
out representation of the type
$\theta_1(x|\tau)=\sum\,C_{mnp}\,g_2^m\,g_3^n\,\eta^p\,x^k$ through
the recurrence of Weierstrass $A_{m,n}$, but $\mathfrak{G}_{m,n}$ is
more effective than  $A_{m,n}$ since polynomials are already grouped
together in $\vartheta$-constants.

Odd derivatives $\theta_1^{\sss(2k+1)}(0|\tau)=
\vartheta_1^{\sss(2k+1)}(\tau)$, i.\,e. expressions in front of
$x^{2k+1}$ in (\ref{t1}), generate polynomials (\ref{theta1}) in
$(\eta,\vartheta)$ which are exactly integrable $k$ times in $\tau$.

\subsection{Functions $\theta_{2,3,4}(x|\tau)$}
Expansions of functions $\thetaAB{\alpha}{\beta}=\pm\theta_{2,3,4}$
of the form
\begin{equation}\label{theta234}
\thetaAB{\alpha}{\beta}(x|\tau) =
\sum_{k=0}^{\infty}\,C_k^{\sss(\alpha,\beta)}(\tau)\!\cdot\!x^{2k}=
\varthetaAB{\alpha}{\beta}- \varthetaAB{\alpha}{\beta} \left\{
2\,\eta+\mfrac{\pi^2}{6}
\Big(\spin{\beta}\,\varthetaAB{\alpha\mbox{\tiny--}1}{0}^4-
\spin{\alpha}\,\varthetaAB{0}{\beta\mbox{\tiny--}1}^4\Big)
\right\}x^2+\cdots
\end{equation}
have the following analytic representation:
\begin{equation}\label{series}
\thetaAB{\alpha}{\beta}(x|\tau) = \mbox{\large$\ds\sum_{\scriptstyle
k=0}^\infty$} \,\frac{(4\pi\ri)^k}{(2k)!}\,
\frac{d^k\varthetaAB{\alpha}{\beta}}{d\tau^k} \!\cdot\!x^{2k}\;.
\end{equation}
The equation (\ref{epsilon}) on functions
$\Xi=(\sigma,\sigmalambda)$ in $(\alpha,\beta)$-representation has
the form
$$
\frac{\partial^2\,\Xi}{\partial
x^2}+\mbox{\large$\widehat{\boldsymbol{\mathfrak D}}$}\, \Xi
+\left\{e_{\sss\gamma\delta}(\vartheta) + \mfrac{\pi^4}{12^2} \big(
\vartheta_{\sss\alpha0}{}^{\!\!\!\!\!\!8}\; +\spin{\alpha+\beta}\,
\vartheta_{\sss\alpha0}{}^{\!\!\!\!\!\!4}\;\;
\vartheta_{\sss0\beta}{}^{\!\!\!\!\!\!4}\;+
\vartheta_{\sss0\beta}{}^{\!\!\!\!\!\!8}\; \big) x^2
\right\}\Xi=0\;,
$$
where $e_{\sss\gamma\delta}(\vartheta)$ correspond  (independently
of representation $(\alpha,\beta)$) to chosen function  $\sigma$ or
$\sigmalambda$:

\begin{equation}\label{ee}
e_{\sss\gamma\delta}(\vartheta)= \mfrac{\pi^2}{12}\! \left(
\spin{\delta}\,\varthetaAB{\gamma\mbox{\tiny--}1}{0}^4-
\spin{\gamma}\,\varthetaAB{0}{\delta\mbox{\tiny--}1}^4\right)\;,
\qquad\qquad \left\{
\begin{array}{ll}
\ds e_{\sss 11}\equiv 0,&
e_{\sss 10}\equiv e_{\sss 1}\\
e_{\sss 00}\equiv e_{\sss 2},& e_{\sss 01}\equiv e_{\sss 3}
\end{array}
\right\}\;.
\end{equation}
Representation is called symmetrical if
$(\gamma,\delta)=(\alpha,\beta)$.

In symmetrical representation the series (\ref{series}) has the
following grouped form:
\begin{equation}\label{t234}
\thetaAB{\alpha}{\beta}(x|\tau)= \varthetaAB{\alpha}{\beta}(\tau)
\,\, \mbox{\large$\ds\sum_{\scriptstyle k=0}^\infty$}
\,(-2)^k\left\{ \sum_{\sss\nu=0}^{\sss k}
\mfrac{(-6)^{-\nu}\,\pi^{2\nu}}{(k-\nu)!\,(2\nu)!} \!\cdot\!
\eta^{k-\nu}\,\mathcal{N}_\nu^{\sss(\alpha,\beta)}(\vartheta)
\right\} x^{2k}
\end{equation}
with the universal integer recurrence $(\mathfrak{G}_{\sss
0,\,0}^{\sss(\alpha,\beta)}=1)$:
$$
\begin{array}{rl}
\ds\! \mathcal{N}_\nu^{\sss(\alpha,\beta)}(\vartheta)\!\!&=\,
\mbox{\large$\ds\sum_{s=0}^{\nu}$}\;
\mathfrak{G}_{\nu-s,\,s}^{\sss(\alpha,\beta)} \cdot
\varthetaAB{0}{\beta\mbox{\tiny--}1}^{4s\mathstrut}\,
\varthetaAB{\alpha\mbox{\tiny--}1}{0}^{4{\sss(}\nu-s{\sss)}}\;,\\ \\
\ds \mathfrak{G}_{m,\,n}^{\sss(\alpha,\beta)}\!\!&=\,
\spin{\alpha}\,(4\,n-8\,m-3)\,
\mathfrak{G}_{m,\,n-1}^{\sss(\alpha,\beta)}-
\spin{\beta}\,(4\,m-8\,n-3)\,
\mathfrak{G}_{m-1,\,n}^{\sss(\alpha,\beta)}-
\\ \\&\ds\phantom{=\,}
-2\,(m+n-1)(2\,m+2\,n-3)
\big(\mathfrak{G}_{m-2,\,n}^{\sss(\alpha,\beta)}+
\spin{\alpha+\beta}\, \mathfrak{G}_{m-1,\,n-1}^{\sss(\alpha,\beta)}+
\mathfrak{G}_{m,\,n-2}^{\sss(\alpha,\beta)}\big)\;.
\end{array}
$$
Symmetry of the recurrence
$\mathfrak{G}_{m,\,n}^{\sss(\alpha,\beta)}$ with respect to
permutation of indices is determined by the properties:
$$
\begin{array}{rl}
\ds \mathfrak{G}_{n,\,m}^{\sss(\alpha,\beta)}=
(-1)^{\sss(m+n)(\alpha+\beta+1)}\,
\mathfrak{G}_{m,\,n}^{\sss(\alpha,\beta)}\;,\\
\mathfrak{G}_{m,\,n}^{\sss(\beta,\alpha)}=
(-1)^{\sss(m+n)(\alpha+\beta)\phantom{+1}}\,
\mathfrak{G}_{m,\,n}^{\sss(\alpha,\beta){}^{\ds\mathstrut}}\;.
\end{array}
$$
This means that there are only two recurrence
$\mathfrak{G}_{m,\,n}^{\sss(\alpha)}$ under $\alpha=\{1,0\}$ and
$\beta=0$:
$$
\begin{array}{rl}
\ds \mathfrak{G}_{m,\,n}^{\sss(\alpha)}\!\!&=\,
\spin{\alpha}\,(4\,n-8\,m-3)\,
\mathfrak{G}_{m,\,n-1}^{\sss(\alpha)}- (4\,m-8\,n-3)\,
\mathfrak{G}_{m-1,\,n}^{\sss(\alpha)}-
\\ \\&\phantom{=\,}
-2\,(m+n-1)(2\,m+2\,n-3) \big(\mathfrak{G}_{m-2,\,n}^{\sss(\alpha)}+
\spin{\alpha}\, \mathfrak{G}_{m-1,\,n-1}^{\sss(\alpha)}+
\mathfrak{G}_{m,\,n-2}^{\sss(\alpha)}\big)
\end{array}
$$
with property $\mathfrak{G}_{n,\,m}^{\sss(\alpha)}=
(-1)^{\sss(m+n)(\alpha+1)}\, \mathfrak{G}_{m,\,n}^{\sss(\alpha)}$:
$$
\mathfrak{G}_{n,\,m}^{\sss(0)}= (-1)^{\sss(m+n)}\,
\mathfrak{G}_{m,\,n}^{\sss(0)}\;,\qquad\qquad
\mathfrak{G}_{n,\,m}^{\sss(1)}= \mathfrak{G}_{m,\,n}^{\sss(1)}\;.
$$
Developments (\ref{t1}) and (\ref{t234}) differ from each other only
in a multiplier and shape of the recurrence. They can be unified
into one (by introducing the parity $\varepsilon$) but the quantity
$\spin{\alpha}$ is left and $(m,n)$-entries of matrices
$\mathfrak{G}^{\sss(\beta,\alpha)}$ differ only in a sign.

Even derivatives $\theta_{\sss\alpha\beta}^{\sss(2k)}(0|\tau)=
\vartheta_{\sss\alpha\beta}^{\sss(2k)}(\tau)$, i.\,e. expressions in
front of $x^{2k}$ in (\ref{t234}), generate polynomials
(\ref{theta234}) in $(\eta,\vartheta)$ which are integrable $k$
times in $\tau$. Their integrability is a consequence of one
dynamical system considered in sect.\,7.

Making use of these expansions one can build other ones about points
$x=\{\pm\frac12,\,\pm\frac\tau2\}$. Described series will be, up to
obvious modifications, changed into one another.

\section{Dynamical systems on $\theta$-functions}

\subsection{Differential equations in $x$}
The five functions
$$
\theta_{1,2,3,4}^{}\quad\mbox{and}\quad\theta_{1\!\!\!}'\equiv
\frac{\partial\theta_1^{}}{\partial x}
$$
satisfy closed ordinary autonomous differential equations:
$$
\left\{
\begin{array}{l}
\ds\frac{\partial\theta_2^{}}{\partial x}=
\frac{\theta_{1\!\!\!}'}{\theta_1^{}}\,\theta_2^{}-
\pi\,\vartheta_2^2\!\cdot\!
\frac{\theta_3^{}\theta_4^{}}{\theta_1^{}}
\\\\\ds
\frac{\partial\theta_3^{}}{\partial x}=
\frac{\theta_{1\!\!\!}'}{\theta_1^{}}\,\theta_3^{}-
\pi\,\vartheta_3^2\!\cdot\!
\frac{\theta_2^{}\theta_4^{}}{\theta_1^{}}
\\\\\ds
\frac{\partial\theta_4^{}}{\partial x}=
\frac{\theta_{1\!\!\!}'}{\theta_1^{}}\,\theta_4^{}-
\pi\,\vartheta_4^2\!\cdot\!
\frac{\theta_2^{}\theta_3^{}}{\theta_1^{}}
\end{array}\right.
\qquad\qquad
\begin{array}{l}
\ds
\frac{\partial\theta_{1\!\!\!}'}{\partial x}=
\frac{\theta_{1\!\!\!}'^2}{\theta_1^{}}-\pi^2\vartheta_3^2\,\vartheta_4^2
\!\cdot\!
\frac{\theta_2^2}{\theta_1^{}}-
\Big\{4\eta+\mfrac{\pi^2}{3}\big(\vartheta_3^4+\vartheta_4^4\big) \Big\}
\!\cdot\!\theta_1^{}
\\\\\ds
\frac{\partial\theta_1^{}}{\partial x}=\theta_{1\!\!\!}'\;. \\ \\
\ds\phantom{\frac{\theta_{1\!\!\!}'}{\theta_1^{\mathstrut}}\,\theta_3^{}}
\end{array}
$$
These equations are equivalent to relations between Weierstrassian
functions $(\sigma,\zeta,\wp,\wpp)(x|\tau)$. General form of
differentiation of the functions $\theta_{1,2,3,4}$ is as follows
($\vartheta_1=0$):
$$
\frac{\partial\theta_k}{\partial x}= \frac{\theta_1'}{\theta_1^{}}
\,\theta_k- \pi\,\vartheta_k^2\!\cdot\!
\frac{\theta_\nu\,\theta_\mu}{\theta_1^{}}\;,\qquad\quad
\mbox{where}\quad\nu=\frac{8\,k-28}{3\,k-10},\quad
\mu=\frac{10\,k-28}{3\,k-8},\quad k=(1,2,3,4)\,.
$$
It immediately follows that the following identities are obeyed
\footnote{These fundamental relations appear implicitly in works by
Jacobi but have not got into comprehensive handbook for elliptic
functions compiled by Schwarz on the basis of Weierstrass's
lectures.}
$$
\frac{\theta_n'}{\theta_n}-\frac{\theta_m'}{\theta_m}=
\pi\vartheta_k^2\!\cdot\!\frac{\theta_1^{}\theta_k}{\theta_n\theta_m}
\,\mbox{\small sign$(n-m)$}\;,\qquad\ \mbox{where}\quad (k\ne n,\; n
\ne m, \; k \ne m, \quad k,n,m=2,3,4)\;,
$$
and also well-known polynomial identities:
\begin{equation}\label{ident}
\theta_1^4+\theta_3^4=\theta_2^4+\theta_4^4,\quad
\ds\vartheta_3^2\,\theta_3^2=\vartheta_2^2\,\theta_2^2+
\vartheta_4^2\,\theta_4^2,\quad  \mbox{\small
sign$(n-m)$}\!\cdot\!\vartheta_k^2\,\theta_1^2=
\vartheta_m^2\,\theta_n^2-\vartheta_n^2\,\theta_m^2\;.
\end{equation}

\subsection{Differential equations in $\tau$}

The functions $\theta_{1,2,3,4}, \theta_{1\!\!\!}'$ satisfy closed
non-autonomous ordinary differential equations in $\tau$:
\begin{equation}\label{Dtau}
\left\{
\begin{array}{l}
\ds \frac{\partial \theta_1^{}}{\partial\tau}=
\mfrac{-\ri}{4\pi}\,\frac{\theta_{1\!\!\!}'^2}{\theta_1^{}}+
\phantom{xxxxxxxxxxx\,}+ \mfrac{\pi
\ri}{4}\,\vartheta_3^2\,\vartheta_4^2\!\cdot\!
\frac{\theta_2^2}{\theta_1^{}}\;\;\;+
\Big\{\mfrac{\ri}{\pi}\eta+\mfrac{\pi \ri}{12}
\big(\vartheta_3^4+\vartheta_4^4\big) \Big\}\!\cdot\!\theta_1^{}\\\\
\ds \frac{\partial \theta_2^{}}{\partial\tau}= \mfrac{-\ri}{4\pi}\!
\left\{\frac{\theta_{1\!\!\!}'}{\theta_1^{}}-
\pi\vartheta_2^2\!\cdot\! \frac{\theta_3^{}\theta_4^{}}
{\theta_1^{}\theta_2^{}}\right\}^{\!\!2}\theta_2^{} + \mfrac{\pi
\ri}{4}\,\vartheta_3^2\,\vartheta_4^2\!\cdot\!
\frac{\theta_1^2}{\theta_2^{}}\;\;\;+
\Big\{\mfrac{\ri}{\pi}\eta+\mfrac{\pi \ri}{12}
\big(\vartheta_3^4+\vartheta_4^4\big) \Big\}\!\cdot\!\theta_2^{}\\\\
\ds \frac{\partial \theta_3^{}}{\partial\tau}= \mfrac{-\ri}{4\pi}\,
\frac{\theta_{1\!\!\!}'^2}{\theta_1^2}\,\theta_3^{}
+\mfrac{\ri}{2}\vartheta_3^2\!\cdot\!
\theta_2^{}\theta_4^{}\frac{\theta_{1\!\!\!}'} {\theta_1^2} -
\mfrac{\pi\ri}{4}\,\vartheta_2^2\,\vartheta_3^2\!\cdot\!
\frac{\theta_4^2}{\theta_1^2}\,\theta_3^{}+
\Big\{\mfrac{\ri}{\pi}\eta+\mfrac{\pi \ri}{12}
\big(\vartheta_3^4+\vartheta_4^4\big) \Big\}\!\cdot\!\theta_3^{}\\\\
\ds \frac{\partial \theta_4^{}}{\partial\tau}=
\mfrac{-\ri}{4\pi}\,\frac{\theta_{1\!\!\!}'^2}{\theta_1^2}\,\theta_4^{}
+\mfrac{\ri}{2}\vartheta_4^2\!\cdot\!
\theta_2^{}\theta_3^{}\frac{\theta_{1\!\!\!}'} {\theta_1^2}-
\mfrac{\pi\ri}{4}\,\vartheta_2^2\,\vartheta_4^2\!\cdot\!
\frac{\theta_3^2}{\theta_1^2}\,\theta_4^{}+
\Big\{\mfrac{\ri}{\pi}\eta+\mfrac{\pi \ri}{12}
\big(\vartheta_3^4+\vartheta_4^4\big) \Big\}\!\cdot\!\theta_4^{}\\\\
\ds \frac{\partial \theta_{1\!\!\!}'}{\partial\tau}=
\mfrac{-\ri}{4\pi}\,\frac{\theta_{1\!\!\!}'^3}{\theta_1^2}
+3\!\left\{ \mfrac{\pi\ri}{4}\vartheta_3^2\,\vartheta_4^2\!\cdot\!
\frac{\theta_2^2}{\theta_1^2} +\mfrac{\ri}{\pi}\eta+
\mfrac{\pi\ri}{12}\big(\vartheta_3^4+\vartheta_4^4\big)
\right\}\theta_{1\!\!\!}'
-\mfrac{\pi^2}{2}\ri\,\vartheta_2^2\,\vartheta_3^2\,\vartheta_4^2\!\cdot\!
\frac{\theta_2^{}\theta_3^{}\theta_4^{}}{\theta_1^2}
\end{array}\right..
\end{equation}
General form of $\tau$-differentiation of the functions
$\theta_{1,2,3,4}$ is as follows:
$$
\begin{array}{rl}
\ds \frac{\partial\theta_k}{\partial \tau}\!\!&=\ds\;
\mfrac{-\ri}{4\pi}\, \mfrac{\theta_1'{}^2}{\theta_1^2} \,\theta_k+
\mfrac{\ri}{2}\, \vartheta_k^2\!\cdot\! \theta_\nu\,
\theta_\mu\mfrac{\theta_1'{}}{\theta_1^2}+\mfrac{\pi\ri}{4}\left\{
\vartheta_3^2\,\vartheta_4^2\!\cdot\!\theta_2^2-
\vartheta_k^2\,\vartheta_\nu^2\,\vartheta_\mu^2\!\cdot\! \left(
\mfrac{\theta_\nu^2}{\vartheta_\nu^2}+
\mfrac{\theta_\mu^2}{\vartheta_\mu^2}\right)\right\}\theta_k
+{}\\ \\
&\ds\phantom{=\;} {}+\left\{
\mfrac{\ri}{\pi}\,\eta+\mfrac{\pi\ri}{12} \big(\vartheta_3^4+
\vartheta_4^4\big)\right\}\!\cdot\! \theta_k\,,\qquad\quad
\mbox{where}\quad\nu=\mfrac{8\,k-28}{3\,k-10},\quad
\mu=\mfrac{10\,k-28}{3\,k-8},\quad k=(1,2,3,4)\,.
\end{array}
$$

\section{Differential equations on $\vartheta,\eta$-constants}

\noindent Differential closure is provided by the constants
$\vartheta_{2,3,4}$ and $\eta$:
\begin{equation}\label{var}
\begin{array}{l}
\ds \frac{1}{\vartheta_{2\!}^{}}\,\frac{d\vartheta_{2\!}^{}}{d\tau}
_{\ds\mathstrut}= \mfrac{\ri}{\pi}\,\eta+\mfrac{\pi\,\ri}{12}\,
\big(\vartheta_{3\!}^4+\vartheta_{4\!}^4 \big)\,, \qquad \ds
\frac{1}{\vartheta_{4\!}^{}}\,\frac{d\vartheta_{4\!}^{}}{d\tau}
^{\ds\mathstrut}= \mfrac{\ri}{\pi}\,\eta-\mfrac{\pi\,\ri}{12}\,
\big(\vartheta_{2\!}^4+\vartheta_{3\!}^4 \big)\,,
\\
\ds \frac{1}{\vartheta_{3\!}^{}}\,\frac{d\vartheta_{3\!}^{}}{d\tau}
_{\ds\mathstrut}^{\ds\mathstrut}=
\mfrac{\ri}{\pi}\,\eta+\mfrac{\pi\,\ri}{12}\,
\big(\vartheta_{2\!}^4-\vartheta_{4\!}^4 \big)\,,\qquad
\quad\;\;\,\frac{d\eta}{d\tau}=\mfrac{\ri}{\pi}\Big\{2\,\eta^2-\mfrac{\pi^4}{12^2}
\big(\vartheta_{2\!}^8+\vartheta_{3\!}^8+\vartheta_{4\!}^8 \big)
\Big\}.
\\
\end{array}
\end{equation}
Well known differential equations on logarithms of ratios
$\vartheta_{2\!}^{}:\vartheta_{3\!}^{}:\vartheta_{4\!}^{}$ and also
known dynamical system of Halphen and its numerous varieties are
widely encountered in the modern literature. These are generated by
the system of equations
$$
\frac{d g_2^{}}{d\tau} = \frac{\ri}{\pi}
\Big(8\,g_2^{}\,\eta-12\,g_3^{}\Big)\,,\qquad \ds \frac{d
g_3^{}}{d\tau} = \frac{\ri}{\pi}
\Big(12\,g_3^{}\,\eta-\mfrac23\,g_2^2\Big)\,,\qquad
\frac{d\eta}{d\tau}=
\frac{\ri}{\pi}\Big(2\,\eta^2-\mfrac16\,g_2^{}\Big)\,,
$$
which, in implicit form, was written out by Weierstrass and Halphen
used it \cite[{\bf I:} p.\,331]{halphen} in order to get his famous
symmetrical version of the system. Ramanujan obtained its equivalent
for some number-theoretic  $q$-series. See for example
\cite[\S\,1]{zudilin} and additional information in this work. It is
lesser known, to all appearances, that another version of the system
was written out by Jacobi in connection with power series of
$\theta$-functions. These results were published by Borchardt on the
basis of papers kept after  Jacobi  \cite[{\bf II}:
383--398]{jacobi}.

Jacobi obtained a dynamical system on four functions $(A,B,a,b)$
which are, in our notation, rational functions of the
$\vartheta,\eta$-constants
$$
A=\vartheta_3^2\,,\quad
B=\frac{4}{\pi^2}\,\frac{\eta}{\vartheta_3^2}
+\frac13\,\frac{\vartheta_2^4-\vartheta_4^4}{\vartheta_3^2}\,,\quad
a=4\!\left(1-2\,\frac{\vartheta_2^2}{\vartheta_3^2}\right)\,,\quad
b=2\,\frac{\vartheta_2^2\vartheta_4^2}{\vartheta_3^4}
$$
and showed that the dynamical system can be useful for obtaining
power series of $\theta$-functions. He also presented two
(canonical) transformations of the variables  $(A,B,a,b)$ retaining
the shape of the equations. Jacobi noticed that the developments
become simple and have a recursive form when extracted an
exponential multiplier $\re^{-\frac12 ABx^2}$. This is, in fact,
Halphen's recurrence (\ref{halphen}) for the $\sigma$-function of
Weierstrass. Transition between variables in mentioned systems is
not one-to-one but always algebraic. This is because all of them are
consequences of the equations (\ref{var}) since the variables,
appeared in these systems, are rationally expressible through the
$\vartheta, \eta$-constants. A related dynamical system arose  in a
work of Jacobi earlier, when he derived his known equation on
$\vartheta$-constants \cite[{\bf II:} p.\,176]{jacobi}
$$
\big(\vartheta^2\,\vartheta_{\tau\tau\tau}-15\,
\vartheta\,\vartheta_\tau \vartheta_{\tau\tau}+
30\,{\vartheta_\tau}^{\!\!\!3}\big)^2
+32\,\big(\vartheta\,\vartheta_{\tau\tau}-3\,{\vartheta_\tau}^{\!\!\!2}
\big)^3= -\pi^2\,\vartheta^{10}\big(\vartheta\,\vartheta_{\tau\tau}-
3\,{\vartheta_\tau}^{\!\!\!2} \big)^2\,.
$$
Logarithmic derivatives of $\vartheta$-constants, in turn, satisfy
compact differential equation
$$
\big(f_\tau^{}-2\,f^2 \big)f_{\tau\tau\tau}^{} -f_{\tau\tau}^2
+16\,f^3f_{\tau\tau}^{}+4\,f_\tau^2\big(f_\tau^{}-6\,f^2 \big)=0,
\qquad f=\frac{d}{d\tau}\ln\vartheta_{2,3,4}^{}(\tau)
$$
and Dedekind's function $\Lambda=\ln\ded(\tau)$ satisfies the
differential equation of the 3-rd order
$$
\big\{\Lambda_{\tau\tau\tau}-12\,\Lambda_{\tau\tau}\,\Lambda_{\tau}
+16\,{\Lambda_\tau}^{\!\!\!\!3} \big\}^2
+32\,\big\{\Lambda_{\tau\tau}-2\,{\Lambda_\tau}^{\!\!\!\!2}\big\}^3=
\mfrac{4}{27}\,\pi^6\,\ded^{24}\,.
$$

\section{Modular transformations of functions $\theta_1$ and
$\ded$}

\noindent The modular transformation is necessary for computations
of $\theta,\ded$-functions. General
$\mathrm{PSL}_2(\mathbb{Z})$-transformation for the function
$\theta_{1\!}^{}$ is closed in itself (in contrast to
$\theta_{2,3,4}^{}$):
$$
\left\{
\begin{array}{l}
\ds \theta_{1\!}^{}\Big(\mfrac{x}{c\,\tau+d}\Big|
\mfrac{a\,\tau+b}{c\,\tau+d}\Big)= \mbox{\Large
$\boldsymbol\varepsilon$}_\theta^{}(a,c,d)
\!\cdot\!\sqrt{c\,\tau+d\,}\,\, \mbox{\large
$\re$}^{\frac{\scriptstyle\pi\ri \,c\, x^2} {\scriptstyle
c\,\tau+d}}_{\mathstrut}\, \theta_{1\!}^{}(x|\tau),\qquad
\mbox{\Large $\boldsymbol\varepsilon$}_\theta^{}(a,c,d)=\sqrt[8]{1}\\\\
\ds
\mbox{\Large $\boldsymbol\varepsilon$}_\theta^{}(a,c,d)=
\mbox{\Large $\re$}_{\mathstrut}^{3\pi\ri
\left\{ \phantom{\sum\limits_k^{m^2}}\right.
\!\!\!\!\!\!\!\!\!\frac{a-d}{12\,c}
-\frac{d}{6}(2c-3)-\frac14-\frac{c-1}{4}\mbox{\scriptsize sign}
(\mbox{\tiny--}d)+\frac1c\,
{\textstyle\sum\limits_
{\scriptscriptstyle k=1}^{\scriptscriptstyle c-1}}\,k
\left[\frac{d}{c}k\right]\!\!\!\!\!\!\!\!\!
\left.\phantom{\sum\limits_k^{m^2}}\right\}},
\qquad
\mbox{\small under normalization}\quad c>0\\\\
\theta_{1\!}^{}(x|\tau+N)=
\re_{\mathstrut}^{\frac{\pi\ri}{4}\mbox{\tiny $N$}}
\!\!\cdot\!\theta_{1\!}^{}(x|\tau)
\end{array}
\right.
$$
When discussing the modular properties of $\theta$-functions one is
usually pointed out ``where \mbox{\Large
$\boldsymbol\varepsilon_\theta^{}$} is an eighth root of unity'' but
complete algorithm of its computation leaves aside. Meantime, as it
follows from the shape of  $\theta$-series, it is seen that the
series, which can be transformed into  a ``hyper-convergent'' form,
can turn out to be incomputable. It is known that the multiplier
 \mbox{\Large
$\boldsymbol\varepsilon_\theta^{}(a,c,d)$} can be expressed via
Jacobi's symbol $\big(\frac ab\big)$ (Hermite). It has independent
rules of computations. These formulas are displayed  here for
completeness of computations and contain some simplifications of
Dedekind's sums (see \cite{apostol} about them):
$$
\left\{
\begin{array}{rl}
\ds \ded\!\left(\mfrac{a\,\tau+b}{c\,\tau+d}\right)\!\!&=\;
\mbox{\Large $\re$}_{\mathstrut}^{\pi\ri \left\{
\phantom{\sum\limits_k^{m^2}}\right.
\!\!\!\!\!\!\!\!\!\frac{a-d}{12\,c}
-\frac{d}{6}(2c-3)-\frac14-\frac{c-1}{4}\mbox{\scriptsize sign}
(\mbox{\tiny--}d)+\frac1c\, {\textstyle\sum\limits_
{\scriptscriptstyle k=1}^{\scriptscriptstyle c-1}}\,k
\left[\frac{d}{c}k\right]\!\!\!\!\!\!\!\!\!
\left.\phantom{\sum\limits_k^{m^2}}\right\}}\,
\sqrt{c\,\tau+d\,}\,\,\ded(\tau)\\ \\
\ded(\tau+N)\!\!&=\;\mbox{\large $\re$}^{\frac{\pi\ri}{12}
\mbox{\tiny$N$}}_{\ds\mathstrut}\,\ded(\tau)
\end{array}\right.\;.
$$

\section{$\theta$-functions with characteristics}

\noindent The preceding  results, making use of notation with
characteristics, allow us to unify formulas and can be a subject for
further generalizations to higher genera.

\subsection{$(\alpha,\beta)$-representations}

Any object, symmetrical in $\vartheta$-constants, can be rewritten
in $(\alpha,\beta)$-representation. For example representation of
branch points through the $\vartheta$-constants (\ref{ee}).
$\vartheta_{\sss\alpha\beta}$-representations for quantities
$g_{2,3}^{}$ have the following form:
$$
\begin{array}{l}
\ds g_2^{}(\tau)=\phantom{1} \mfrac{\pi^4}{12} \left\{
\vartheta_{\sss\alpha0}{}^{\!\!\!\!\!\!8}\;+ \spin{\alpha+\beta}\,
\vartheta_{\sss\alpha0}{}^{\!\!\!\!\!\!4}\;\;
\vartheta_{\sss0\beta}{}^{\!\!\!\!\!\!4}\;+
\vartheta_{\sss0\beta}{}^{\!\!\!\!\!\!8}\;
\right\}\;, \qquad\qquad (\alpha,\beta)\ne(0,0)\\ \\
\ds g_3^{}(\tau)=\frac{\pi^6}{432}\Big\{
2\,\spin{\beta}\,\vartheta_{\sss\alpha0}{}^{\!\!\!\!\!\!12} -
3\,\vartheta_{\sss\alpha0}^4\, \vartheta_{\sss0\beta}^4
\big(\spin{\beta}\,\vartheta_{\sss0\beta}^4-\spin{\alpha}\,
\vartheta_{\sss\alpha0}^4 \big)
-2\,\spin{\alpha}\,\vartheta_{\sss0\beta}{}^{\!\!\!\!\!\!12}
 \Big\}\;.
\end{array}
$$
Identity $\vartheta_3^4=\vartheta_2^4+\vartheta_4^4$ and Jacobi's
formula $\vartheta_1'=2\pi\ded^3$, in
$(\alpha,\beta)$-representation, have the form
\begin{equation}\label{jacobi}
\begin{array}{l}
\;\;\,\varthetaAB{\alpha}{\beta}^4= \left(
\spin{\beta}\,\varthetaAB{\alpha\mbox{\tiny--}1}{0}^4+
\spin{\alpha}\,\varthetaAB{0}{\beta\mbox{\tiny--}1}^4\right)
\mfrac{\spin{\alpha\beta}+1}{2}\\ \\
\vartheta_{\sss\!\alpha\beta}{}{\!\!\!\!'}\,\,(\tau)=
\ri^\beta\,(1-\spin{\alpha\beta})\!\cdot\!
\pi\,\ri\,\ded^3(\tau)\;,\qquad \mbox{where\ \ }
\vartheta_{\sss\!\alpha\beta}{}{\!\!\!\!'}\,\,(\tau)\equiv
\theta_{\sss\!\alpha\beta}{}{\!\!\!\!'}\,\,(0|\tau)\;.
\end{array}
\end{equation}
In general case, when the $\vartheta'$-constant is a value of the
derivative of $\theta$-function at some $\frac12$-period, the
formula of Jacobi (\ref{jacobi}) is generalized to the following
expression:
\begin{equation}\label{jac}
\theta_{\sss\!\alpha\beta}{}\!\!\!\!'\,\, \Big(\mfrac n2+\mfrac
m2\,\tau \Big|\tau\Big)= (-\ri)^{\sss
m(\beta+n)}_{\mathstrut}\!\cdot\! \pi\,\ri\,\re^{\!-\frac
{\pi\ri}{4}m^2\tau}_{\mathstrut} \left\{
\ri^{\beta+n}_{\mathstrut}\,\big(1-\spin{\alpha+m}^{\beta+n}  \big)
\!\cdot\!\ded^3-m\,
\varthetaAB{\alpha\mbox{\tiny+}m}{\beta\mbox{\tiny+}n} \right\}\;.
\end{equation}

Algebraic and differential closeness of Jacobi's functions with
integral characteristics lead to the fact that all the derivatives
$\theta_{\sss\!\alpha\beta}{}\!\!\!\!'\,\,(x|\tau)$, with arguments
shifted by arbitrary $\frac12$-periods, are expressible in terms of
the function $\theta_{1\!}'(x|\tau)$ and functions
$\theta_{\sss\!\alpha\beta\!}(x|\tau)$:
$$
\begin{array}{rl}
\theta_{\sss\!\alpha\beta}{}\!\!\!\!'\,\, \Big(x+\mfrac n2+\mfrac
m2\,\tau \Big|\tau\Big)\!\!&= (-\ri)^{\sss
m(\beta+n)}_{\mathstrut}\!\cdot\! \re^{\!-\pi\ri m\left(x+\frac
m4\tau \right)}_{\mathstrut} \left\{ \Big(
\mfrac{\theta_{1\!}'(x|\tau)}{\theta_{1\!}^{}(x|\tau)}-\pi\ri\,m
\Big)\,\thetaAB{\alpha\mbox{\tiny+}m}{\beta\mbox{\tiny+}n}(x|\tau)-{}
\right.\\ \\& \ds\phantom{=}\left.
{}-\spin{\alpha+m}^{\left[\frac{\beta+n}{2}\right]}_{\mathstrut}
\!\cdot\!\pi\,
\varthetaAB{\alpha\mbox{\tiny+}m}{\beta\mbox{\tiny+}n}^2
\!\cdot\!\mfrac{\thetaAB{1\mbox{\tiny--}\alpha\mbox{\tiny--}m}
{0}(x|\tau)\,
\varthetaAB{0}{1\mbox{\tiny--}\beta\mbox{\tiny--}n}(x|\tau)}
{\theta_{1\!}^{}(x|\tau)} \right\}\;.
\end{array}
$$
Under $x=0$ need to use the formula (\ref{jac}).

\subsection{Modular transformations}
The modular transformations served as the key for Hermite when he
obtained his celebrated solution of the quintic equation
$x^5-x-A=0$\footnote{Solution of this equation is expressed, in
fact, in terms of $\vartheta$-constants and was obtained by Hermite
in 1858 (see last formula on p.\,10 in \cite{hermite}). Curiously
that it contains an erroneous sign which is repeated, to the best of
our knowledge, everywhere the solution is reproduced.}. Hermite
wrote out the multiplier $\mbox{\Large
$\boldsymbol\varepsilon$}_\theta^{}(a,c,d)$ in a form of Gaussian
sum of exponents. We propose the exponent of a sum what is simpler.
General modular transformation $\tau \to \frac{a\,\tau+b}{c\,\tau+d}
$ for the functions $\theta$ has the form
$$
\left\{
\begin{array}{l}
\ds \theta \raisebox{0.04em}{\mbox{$ \scalebox{0.5}[0.9]{\big[}
\!\!\raisebox{0.1em} {\scalebox{0.8}{\mbox{\tiny$
\begin{array}{c}\alpha'\mbox{\tiny--}1\\\beta'\mbox{\tiny--}1
\end{array}$}}}
\!\! \scalebox{0.5}[0.9]{\big]}$}}
\Big(\mfrac{x}{c\,\tau+d}\Big|\mfrac{a\,\tau+b}{c\,\tau+d}\Big)=
\mbox{\Large $\boldsymbol\varepsilon$}_\theta^{}(a,c,d)\,
\mbox{\large$\re$}_{\mathstrut}^ {\frac{\pi\ri}{4}\left\{
2\alpha(bc\,\beta-d+1)-c\beta(a\beta-2)-db\,\alpha^2 \right\}
}\!\cdot\! \sqrt{c\,\tau+d\,}\,\,
\mbox{\large$\re$}^{\frac{\pi\ri\,cx^2}{c\,\tau+d}}_{\mathstrut}\,
\theta \raisebox{0.04em}{\mbox{$ \scalebox{0.5}[0.9]{\big[}
\!\!\raisebox{0.1em} {\scalebox{0.8}{\mbox{\tiny$
\begin{array}{c}\alpha\mbox{\tiny--}1\\\beta\mbox{\tiny--}1
\end{array}$}}}
\!\! \scalebox{0.5}[0.9]{\big]}$}}
(x|\tau)\\\\
\ds \mbox{\Large $\boldsymbol\varepsilon$}_\theta^{}(a,c,d)=
\mbox{\Large $\re$}_{\mathstrut}^{3\pi\ri \left\{
\phantom{\sum\limits_k^{m^2}}\right.
\!\!\!\!\!\!\!\!\!\frac{a-d}{12\,c}
-\frac{d}{6}(2c-3)-\frac14-\frac{c-1}{4}\mbox{\scriptsize sign}
(\mbox{\tiny--}d)+\frac1c\, {\textstyle\sum\limits_
{\scriptscriptstyle k=1}^{\scriptscriptstyle c-1}}\,k
\left[\frac{d}{c}k\right]\!\!\!\!\!\!\!\!\!
\left.\phantom{\sum\limits_k^{m^2}}\right\}}, \qquad
\mbox{\small under normalization}\quad c>0\\\\
\theta \raisebox{0.04em}{\mbox{$ \scalebox{0.5}[0.9]{\big[}
\!\!\raisebox{0.1em} {\scalebox{0.8}{\mbox{\tiny$
\begin{array}{c}\alpha\mbox{\tiny--}1\\\beta\end{array}$}}}
\!\! \scalebox{0.5}[0.9]{\big]}$}} (x|\tau+N)=
\re_{\mathstrut}^{\mbox{\tiny--}\frac{\pi\ri}{4} \mbox{\tiny
$N$}(\alpha^2\mbox{\tiny--}1)} \!\cdot\! \theta
\raisebox{0.04em}{\mbox{$ \scalebox{0.5}[0.9]{\big[}
\!\!\raisebox{0.1em} {\scalebox{0.8}{\mbox{\tiny$
\begin{array}{c}\alpha-1\\\beta\!+\!N\alpha\end{array}$}}}
\!\! \scalebox{0.5}[0.9]{\big]}$}}(x|\tau)
\end{array}
\right.
$$
where, at given modular transformation, characteristics are computed
by the formulas:
$$
\left\{
\begin{array}{l}
\ds
\alpha'=\phantom{-}d\,\alpha-c\,\beta\\ \\
\ds \beta'=-b\,\alpha +a\,\beta
\end{array}\right.\,,
\qquad\qquad \left\{
\begin{array}{l}
\ds
\alpha=a\,\alpha'+c\,\beta'\\ \\
\beta=b\,\alpha' +d\,\beta'
\end{array}\right.\;.
$$
Under $x=0$ the transformation turns into the transformation of the
$\vartheta$-constants
$$
\left\{
\begin{array}{l}
\ds \vartheta \raisebox{0.04em}{\mbox{$ \scalebox{0.5}[0.9]{\big[}
\!\!\raisebox{0.1em} {\scalebox{0.8}{\mbox{\tiny$
\begin{array}{c}\alpha'\mbox{\tiny--}1\\\beta'\mbox{\tiny--}1
\end{array}$}}}
\!\! \scalebox{0.5}[0.9]{\big]}$}}
\Big(\mfrac{a\,\tau+b}{c\,\tau+d}\Big)= \mbox{\Large
$\boldsymbol\varepsilon$}_\theta^{}(a,c,d)\,
\mbox{\large$\re$}_{\mathstrut}^ {\frac{\pi\ri}{4}\left\{
2\alpha(bc\,\beta-d+1)-c\beta(a\beta-2)-db\,\alpha^2 \right\}
}\!\cdot\! \sqrt{c\,\tau+d\,}\, \vartheta \raisebox{0.04em}{\mbox{$
\scalebox{0.5}[0.9]{\big[}  \!\!\raisebox{0.1em}
{\scalebox{0.8}{\mbox{\tiny$
\begin{array}{c}\alpha\mbox{\tiny--}1\\\beta\mbox{\tiny--}1
\end{array}$}}}
\!\! \scalebox{0.5}[0.9]{\big]}$}}
(\tau)\\\\
\ds \mbox{\Large $\boldsymbol\varepsilon$}_\theta^{}(a,c,d)=
\mbox{\Large $\re$}_{\mathstrut}^{3\pi\ri \left\{
\phantom{\sum\limits_k^{m^2}}\right.
\!\!\!\!\!\!\!\!\!\frac{a-d}{12\,c}
-\frac{d}{6}(2c-3)-\frac14-\frac{c-1}{4}\mbox{\scriptsize sign}
(\mbox{\tiny--}d)+\frac1c\, {\textstyle\sum\limits_
{\scriptscriptstyle k=1}^{\scriptscriptstyle c-1}}\,k
\left[\frac{d}{c}k\right]\!\!\!\!\!\!\!\!\!
\left.\phantom{\sum\limits_k^{m^2}}\right\}}, \qquad
\mbox{\small under normalization}\quad c>0\\\\
\vartheta \raisebox{0.04em}{\mbox{$ \scalebox{0.5}[0.9]{\big[}
\!\!\raisebox{0.1em} {\scalebox{0.8}{\mbox{\tiny$
\begin{array}{c}\alpha\mbox{\tiny--}1\\\beta\end{array}$}}}
\!\! \scalebox{0.5}[0.9]{\big]}$}} (\tau+N)=
\re_{\mathstrut}^{\mbox{\tiny--}\frac{\pi\ri}{4} \mbox{\tiny
$N$}(\alpha^2\mbox{\tiny--}1)} \!\cdot\!
\varthetaAB{\alpha-1}{\beta\!+\!N\alpha}(\tau)
\end{array}
\right.
$$
If characteristics $(\alpha,\beta)=\{0,1\}\;\mbox{mod}\;2$ then the
formulas are closed: $(\alpha',\beta')=\{0,1\}\;\mbox{mod}\;2$.
Ratio of any two $\vartheta,\theta$-functions contains no the
multiplier \mbox{\Large
$\boldsymbol\varepsilon_\theta^{}$}$(a,c,d)$. Hermite used this fact
to build his functions $\varphi,\psi,\chi(\tau)$ and transformations
between them.

\subsection{Multiplication theorems}

\noindent Let $n$ be arbitrary complex number. Then the functions
$\theta$ satisfy the complex multiplication theorems determined by
the following recurrences :
$$
\left\{
\begin{array}{l}
\ds \theta_1^{}(nx)=\frac{\theta_3^2(n_1^{} x)\,\theta_2^2(x)-
\theta_2^2(n_1^{}x)\,\theta_3^2(x)} {\vartheta_4^2 \cdot
\theta_1^{}\big((n-2)x\big)}\,,
\qquad\theta_1^{}(2x)=2\,\theta_1^{}(x)\,
\frac{\theta_2^{}(x)\,\theta_3^{}(x)\,\theta_4^{}(x)}
{\vartheta_2^{}\vartheta_3^{}\vartheta_4^{}}\\\\
\ds \thetaAB{\alpha+1}{\beta+1}(nx)=-
\frac{\langle\beta\rangle\,\thetaAB{\alpha}{0}^2(n_1^{}x)\,
\thetaAB{\alpha}{0}^2(x)+\langle\alpha\rangle\,
\thetaAB{0}{\beta}^2(n_1^{}x)\,\thetaAB{0}{\beta}^2(x)}
{\varthetaAB{\alpha+1}{\beta+1}^2 \cdot
\thetaAB{\alpha+1}{\beta+1}\big((n-2)x\big)}\,, \qquad n_1^{}\equiv
n-1
\end{array}\right..
$$
This is not the only representation. Multiplication for $\theta_1'$
is derived by taking a derivative. If $n$ is integer then the
formulas are closed. In particular, multiplications for the
functions $\theta_{2,3,4}$ are closed (doublings were written out by
Jacobi):
$$
\left\{
\begin{array}{l}
\ds \theta_2^{}(nx)=\frac{\theta_3^2(n_1^{}x)\,\theta_3^2(x)-
\theta_4^2(n_1^{}x)\,\theta_4^2(x)} {\vartheta_2^2 \cdot
\theta_2^{}\big((n-2)x\big)}\\
\ds
\theta_3^{}(nx)=\frac{\theta_2^2(n_1^{\ds\mathstrut}x)\,\theta_2^2(x)+
\theta_4^2(n_1^{}x)\,\theta_4^2(x)} {\vartheta_3^2 \cdot
\theta_3^{}\big((n-2)x\big)}\\
\ds
\theta_4^{}(nx)=\frac{\theta_3^2(n_1^{\ds\mathstrut}x)\,\theta_3^2(x)-
\theta_2^2(n_1^{}x)\,\theta_2^2(x)} {\vartheta_4^2 \cdot
\theta_4^{}\big((n-2)x\big)}
\end{array}\right.\,.
$$
There exist the polynomial non-recursive multiplications but these
are multiplications not over the constants
$\vartheta_{\sss\alpha\beta}(\tau)$. This follows from
Weierstrassian  formulas for the functions
$\sigma,\sigmalambda(nx)$.

\subsection{Differential equations}

\subsubsection{$\theta$-functions} $\theta$-functions with arbitrary
integer characteristics  $(\alpha,\beta)$, as functions of the two
variables $x$ and $\tau$, satisfy splitted and closed systems of
ordinary differential equations over the differential field of
$\vartheta,\eta$-constants $\mathbb{C}_\partial(\vartheta,\eta)$:
\begin{equation}\label{x}
\left\{
\begin{array}{rl}
\ds \frac{\partial\thetaAB{\alpha}{\beta}}{\partial x}\!\!&=\ds\;
\frac{\theta_1'}{\theta_1^{}} \,\thetaAB{\alpha}{\beta}-
\spin{\alpha}^{\left[\frac{\beta}{2}\right]}_{\mathstrut}\,
\pi\,\varthetaAB{\alpha}{\beta}^2\!\cdot\!
\frac{\thetaAB{1\mbox{\tiny--}\alpha}{0}\,
\thetaAB{0}{1\mbox{\tiny--}\beta}}{\theta_1^{}}\\ \\
\ds \frac{\partial\theta_{1\!\!\!}'}{\partial x}\!\!&=\ds\;
\frac{\theta_{1\!\!\!}'^2}{\thetaAB{1}{1}}-\pi^2
\varthetaAB{0}{0}^2\,\varthetaAB{0}{1}^2 \!\cdot\!
\frac{\thetaAB{1}{0}^2}{\thetaAB{1}{1}}-
\Big\{4\eta+\mfrac{\pi^2}{3}\big(\varthetaAB{0}{0}^4+
\varthetaAB{0}{1}^4\big) \Big\} \!\cdot\!\thetaAB{1}{1}
\end{array}\right.\,,\qquad\qquad\quad\;\;\;
\end{equation}

\begin{equation}\label{tau}
\left\{
\begin{array}{rl}
\ds \frac{\partial\thetaAB{\alpha}{\beta}}{\partial \tau}\!\!&=\ds\;
\mfrac{-\ri}{4\pi}\, \mfrac{\theta_1'{}^2}{\thetaAB{1}{1}^2}
\,\thetaAB{\alpha}{\beta}+
\spin{\alpha}^{\left[\frac{\beta}{2}\right]}_{\mathstrut}\,
\mfrac{\ri}{2}\, \varthetaAB{\alpha}{\beta}^2\!\cdot\!
\mfrac{\thetaAB{1\mbox{\tiny--}\alpha}{0}\,
\thetaAB{0}{1\mbox{\tiny--}\beta}}{\thetaAB{1}{1}^2}\,
\theta_1'{}+{}\\ \\
&\ds\phantom{=\;} {}+\mfrac{\pi\ri}{4}\,
\varthetaAB{1}{0}^2\,\varthetaAB{0}{0}^2\,\varthetaAB{0}{1}^2\!\cdot\!
\left\{ \mfrac{\thetaAB{1}{0}^2}{\varthetaAB{1}{0}^2}-
\mfrac12\,(1+\spin{\alpha\beta})\!\cdot\!\! \left(
\mfrac{\thetaAB{1\mbox{\tiny--}\alpha}{0}^2}
{\varthetaAB{1\mbox{\tiny--}\alpha}{0}^2}+
\mfrac{\thetaAB{0}{1\mbox{\tiny--}\beta}^2}
{\varthetaAB{0}{1\mbox{\tiny--}\beta}^2}
\right)\right\}\thetaAB{\alpha}{\beta}
+{}\\ \\
&\ds\phantom{=\;} {}+\left\{
\mfrac{\ri}{\pi}\,\eta+\mfrac{\pi\ri}{12} \big(\varthetaAB{0}{0}^4+
\varthetaAB{0}{1}^4\big)\right\}\!\cdot\!
\thetaAB{\alpha}{\beta}\\ \\
\ds \frac{\partial \theta_{1\!\!\!}'}{\partial\tau}\!\!&=\ds\;
\mfrac{-\ri}{4\pi}\,\frac{\theta_{1\!\!\!}'^3}{\thetaAB{1}{1}^2}
+3\left\{
\mfrac{\pi\ri}{4}\varthetaAB{0}{0}^2\,\varthetaAB{1}{0}^2\!\cdot\!
\frac{\thetaAB{0}{1}^2}{\thetaAB{1}{1}^2} +\mfrac{\ri}{\pi}\eta+
\mfrac{\pi\ri}{12}\big(\varthetaAB{0}{0}^4+\varthetaAB{0}{1}^4\big)
\right\}\theta_{1\!\!\!}' -{}\\ \\
&\ds \phantom{=}{}-\ri\,\mfrac{\pi^2}{2}\,
\varthetaAB{1}{0}^2\,\varthetaAB{0}{0}^2\,\varthetaAB{0}{1}^2
\!\cdot\! \frac{\thetaAB{1}{0}^{}\thetaAB{0}{0}^{}\thetaAB{0}{1}^{}}
{\thetaAB{1}{1}^2}
\end{array}\right.\,.
\end{equation}
These differentiations are not quite symmetrical. In all likelihood,
the separability of variables in differential properties of the
$\theta$'s is not accidental. Especially in those cases, where exact
solutions of integrable nonlinear partial differential equations are
expressible through the $\theta$-functions. All such solutions,
essentially non-stationary and multi-phase, are consequences of
integrability of the only system of ordinary differential equations
on the five Jacobi's functions $\theta,\theta_{1\!\!\!}'$. They
satisfy the heat equation but it is a partial differential equation.

General solution of dynamical system (\ref{x}) is given by the
formulas
$$
\theta_{\sss\alpha\beta}=a\cdot\theta_{\sss\alpha\beta} \big(
x+b|\tau\big)\,\re^{cx}_{}\,,\quad
\theta_1'=a\cdot\big\{\theta_1'\!( x+b|\tau\big)-c\,\theta_{11}\!(
x+b|\tau)\big\}\,\re^{cx}_{}
$$
with arbitrary functions $a,b,c(\tau)$. The system (\ref{tau}) has a
solution of the form
$$
\theta_{\sss\alpha\beta}=a\cdot\theta_{\sss\alpha\beta}( b|\tau)\,,
\quad \theta_1'=-a'\cdot \theta_{11}\!\big(
b|\tau\big)+a\cdot\theta_1'\!( b|\tau)\,,
$$
with arbitrary functions $a,b(x)$. All the solutions should be
supplemented by the relations (\ref{ident}).

\subsubsection{$\vartheta$-constants}
$(\alpha,\beta)$-representation for the system of equations
(\ref{var}) has the following form
\begin{equation}\label{last}
\left\{
\begin{array}{rl}
\ds
\frac{d\varthetaAB{\alpha}{\beta}}{d\tau}\!\!\!&=\,
\varthetaAB{\alpha}{\beta}
\left\{
\mfrac{\ri}{\pi}\,\eta+\mfrac{\pi\ri}{12}
\Big(\spin{\beta}\,\varthetaAB{1\mbox{\tiny--}\alpha}{0}^4-
\spin{\alpha}\,\varthetaAB{0}{1\mbox{\tiny--}\beta}^4\Big)
\right\}\\ \\
\ds \frac{d\eta}{d\tau}\!\!\!&=\, \mfrac{\ri}{\pi} \Big\{2\,\eta^2-
\mfrac{\pi^4}{72}\! \left( \varthetaAB{\alpha}{0}^8+
\varthetaAB{0}{\beta}^8+ \spin{\alpha+\beta}\,
\varthetaAB{\alpha}{0}^4\, \varthetaAB{0}{\beta}^4 \right)\!
\Big\}\qquad \leftarrow\quad(\alpha,\beta)\ne (0,0)\;.
\end{array}\right.
\end{equation}
Its general solution, with three arbitrary constants $(a:b:c:d)$, is
given by the formulas
$$
\varthetaAB{\alpha}{\beta}=\sqrt[-2]{c\tau+d\,}\cdot
\varthetaAB{\alpha}{\beta}\!\left(\mfrac{a\tau+b}{c\tau+d}\right),\qquad
\eta=(c\tau+d)^{-2}\cdot\eta\!\left(\mfrac{a\tau+b}{c\tau+d}\right)+
\frac12\frac{\pi\ri c}{c\tau+d}
$$
and the first identity of Jacobi (\ref{jacobi}). Right hand sides in
these formulas are the $\vartheta,\eta$-series. Unless one assumes
that the quantities $\eta,\vartheta$ in equations (\ref{tau}) being
the $\eta(\tau),\vartheta(\tau)$-constants, then  the requirement of
differential closeness of the field of coefficients of the system,
due to heat equation, leads to the equations (\ref{last}). One may
also view the equations (\ref{last}) as the integrability condition
(compatibility) of the two systems of equations (\ref{x}) and
(\ref{tau}).

\section{Conclusion}

\noindent Dynamical systems of Jacobi--Halphen, their
generalizations, known as $\mathit{SU}(2)$-invariant self-dual
Einstein's equations \cite{hitchin}, modular solutions of equations
of Painlev\'e $\mathrm P_{\sss\mathrm{VI}}$ (Picard, Okamoto, Manin,
and others)  etc are also consequences of splittability of the
equations on $\theta$-functions. Thus, rules of differentiations
(\ref{x}--\ref{last}) generate mentioned dynamical systems and their
solutions and the quantities $\theta,\eta,\vartheta$ are the
uniformizing variables for them. Higher derivatives of the
$\theta$-functions or $\vartheta$-constants are again the
$\eta,\vartheta,\theta,\theta_1'$-functions in form of rational
polynomials of them. Using these facts one can find new solutions or
simplify known ones. For example solution of equation $\mathrm
P_{\sss\mathrm{VI}}$ obtained by Hitchin in the work \cite{hitchin}.

Modular transformations of sect.\,9.2 generate automorphisms of the
systems (\ref{x}--\ref{last}). For the $\vartheta,\theta$-functions
these are merely linear transformations. For other dynamical
systems, like Jacobi's equations on $(A,B,a,b)$, this yields
non-obvious transformations of dynamical variables and systems into
themselves. For the modular solutions of equations of the Painlev\'e
type the transformations  become lesser obvious and cease to be
canonical and change  shape of the equations.

In conclusion we note that it is of interest  to generalize
preceding results to  $\Theta$-functions of higher genera and their
constants. In some particular cases this problem might be tested in
the following manner. Let us consider a situation when the Jacobian
variety admits a decomposition of two-dimensional $\Theta$-function
to functions of Jacobi. For example
$$
\Theta \Big(
\begin{smallmatrix}z_1^{}\\
\mathstrut z_2^{}\end{smallmatrix}\Big|
\ds\begin{smallmatrix}\tau&\frac12_{}\\\frac12\,&
\mu\end{smallmatrix} \Big)=\mfrac12\,\big(\theta_3^{}\!(z_1^{}|\tau)
+\theta_4^{}\!(z_1^{}|\tau) \big)\,
\theta_3^{}\!\big(z_2^{}\big|\mu\big) +
\mfrac12\,\big(\theta_3^{}\!(z_1^{}|\tau)
-\theta_4^{}\!(z_1^{}|\tau)
\big)\,\theta_4^{}\!\big(z_2^{}\big|\mu\big)\;.
$$
Corresponding ten
$\scalebox{1.7}[1.0]{${{\vartheta}}$}
\raisebox{0.04em}{\mbox{$ \scalebox{0.5}[0.9]{\big[}
\!\!\raisebox{0.1em} {\scalebox{0.8}{\mbox{\tiny
$\begin{array}{c}\boldsymbol{\alpha}\\\boldsymbol{\beta}\end{array}$}}}
\!\! \scalebox{0.5}[0.9]{\big]}$}}$-constants are expressible
through the quantities
$$\vartheta_{2,3,4}(\tau)\,,\quad\vartheta_{2,3,4}(\mu)\,,\quad
\theta_{1,2,3,4}^{}\!\big({\textstyle\frac14}\big|\tau\big) \,,\quad
\theta_{1,2,3,4}^{}\!\big({\textstyle\frac14}\big|\mu\big)\,.
$$
Making use of the relations
$$
\begin{array}{ll}
\ds 2\,\theta_3^4\!\big({\textstyle\frac14}\big)=
\vartheta_4^{}\vartheta_3^3+ \vartheta_3^{}\vartheta_4^3\,,\qquad &
2\,\theta_2^4\!\big({\textstyle\frac14}\big)=
\vartheta_4^{}\vartheta_3^3- \vartheta_3^{}\vartheta_4^3
\,,\\\\
\ds \phantom{2\,}\theta_4^{}\!\big({\textstyle\frac14}\big)=
\theta_3^{}\!\big({\textstyle\frac14}\big)\,,\qquad& \ds
\phantom{2\,}\theta_1^{}\!\big({\textstyle\frac14}\big)=
\frac12\,\frac{\vartheta_2^2\,\vartheta_3^{}\, \vartheta_4^{}}
{\theta_2^{}\!\big(\frac14 \big)\,\theta_3^2\!\big(\frac14 \big)}
\end{array}
$$
(these are obtained from the multiplication theorems) one may
analyze the system of derivatives of the
\scalebox{1.7}[1.0]{${{\vartheta}}$}-constants with respect to
moduli  $(\tau,\,\mu)$. In particular, to check its splittability in
$\tau$ and $\mu$. This fact could be also checked, due to the heat
equation, by investigating closeness of differentiations of sixteen
functions $\Theta \raisebox{0.04em}{\mbox{$
\scalebox{0.5}[0.9]{\big[} \!\!\raisebox{0.1em}
{\scalebox{0.8}{\mbox{\tiny
$\begin{array}{c}\boldsymbol{\alpha}\\\boldsymbol{\beta}\end{array}$}}}
\!\! \scalebox{0.5}[0.9]{\big]}$}}(z_1,z_2)$ with respect to
arguments $(z_1,z_2)$. As it is seen from this example, the closure
does  exist and it being algebraic at the least. Whether it will be
polynomial and, if no, how to do this?


\thebibliography{9}

\bibitem{apostol} {\sc Apostol, N.M.} {\em Modular functions and
Dirichlet Series in Number Theory\/}. Springer--Verlag (1976).

\bibitem{eilbeck}{\sc Eilbeck, C., Enol'skii, V.Z.}
{\em Bilinear operators and power series for the Weierstrass
$\sigma$-function\/}. Journ.Phys. {\bf A}: Math.Gen. (2000), {\bf
33}, 791--794.

\bibitem{halphen}{\sc Halphen, G.-H.}
{\em Trait\'e des Fonctions
Elliptiques et de Leurs Applications\/}. {\bf I--III}.
Gauthier--Villars: Paris (1886--1891).

\bibitem{hermite}{\sc Hermite, C.} {\em \OE uvres}. {\bf II}.
Gauthier--Villars: Paris (1908).

\bibitem{hitchin}{\sc Hitchin, N.} {\em Twistor spaces,
Einstein metrics and isomonodromic deformations.} J.\,Diff.\,Geom.
(1995), {\bf 42}(1), 30--112.

\bibitem{jacobi}{\sc Jacobi, C.} {\em Gesammelte Werke\/}.
{\bf I, II}. Verlag von G.Reimer: Berlin (1882).

\bibitem{tannery}{\sc Tannery, J., Molk, J.}
{\em Elements de la
theorie des fonctions elliptiques\/}. {\bf I--IV}.
Gauthier--Villars: Paris (1893--1902).

\bibitem{we}{\sc Weierstrass, K.} {\em Mathematische Werke\/}.
{\bf II, V}. Mayer \& M\"uller: Berlin (1894).

\bibitem{zudilin}{\sc Zudilin, V.V.}
{\em Thetanulls and differential equations\/}. {Russian Acad. Sci.
Sb. Math. {\bf 191}(12) (2000), 1827--1871.}
\end{document}